\documentclass[12pt]{article}

\usepackage{amsfonts}
\usepackage{amsmath, amssymb,esint,color}
\usepackage{graphicx}

\begin{document}
\title{Numerical studies on the self-similar collapse of the $\alpha$-patches problem}
\author{Ana M. Mancho\\
Instituto de Ciencias Matem\'{a}ticas, CSIC-UAM-UC3M-UCM, \\ C/ Nicol\'as Cabrera 15, Campus Cantoblanco UAM, \\28049, Madrid, Spain.}
\date{ }
\maketitle

\begin{abstract}

This paper studies the dynamical evolution of the $\alpha-$patches problem expressed in self-similar variables.
A numerical algorithm is proposed and these  equations are numerically explored.  { Several benchmarks of the code are 
discussed throughout the paper.
Exact self-similar solutions are described  and are
found to play a role in separating collapsing from non-collapsing initial data: small perturbations around this solution blow up while
others do not.   
Numerical simulations performed near convergent rescaled profiles,  such as those described by C\'ordoba et al. in \cite{pnas},     indicate the absence of  a stationary graph in the neighborhood of the rescaled profiles  
and suggest a  more complex scenario for blow up.}

\end{abstract}

\section{Introduction}

A classical open problem in mathematical fluid mechanics is whether 3D Euler equations may develop
singularities in finite time \cite{Const}. 
One of the scenarios proposed in the past for the formation of singularities
in Euler equations is the vortex patch problem. A vortex patch consists of a 2D  
simply connected  and  bounded region of constant vorticity
which is a  weak solution of the 2D Euler equation.
The work of Chemin \cite{Che} and
Bertozzi and Constantin \cite{BeCo} rigorously proved the global existence
of regular solutions for this case, and therefore in the context of this equation singularities cannot appear.

The 2D surface quasigeostrophic equations have also attracted a lot of interest
mainly because of their similarity to the 3D Euler equations \cite{CMT}  and its physical relevance as a
model for the formation of temperature fronts in some geophysical contexts  \cite{CMTpof,CMT,CFR,Ped}. They model  a two-dimensional incompressible fluid system,  
and since it is a 2D model it is 
more tractable than the full 3D Euler system. A very actively studied question for this system is the formation of singularities in finite time for smooth initial data (see  \cite{jap,diego,cf,dh,ccw,clstw}).

The $\alpha$-patches problem is a family of contour dynamics equations that links
the vortex patches scenario to the evolution of patches in the surface quasi-geostrophic equation.
This problem  has been analyzed   by C\'ordoba {\it et al} \cite{pnas}. 
In this setting, the time evolution of  a point, $\vec{x}$, in any of the contours $C_k$ that surrounds a patch with constant $\theta_k$ is described by the equation,
\begin{eqnarray}
\frac{d\overrightarrow{x}(t)}{dt}=\sum_{k=1}^{P} \frac{\theta_k}{2\pi}  \ointctrclockwise_{C_k(t)}\frac{d\overrightarrow{x}_k}{|\overrightarrow{x}(t) -
\overrightarrow{x}_k(t)|^{\alpha}}, \,\,\,\,\,\, \overrightarrow{x} \in \mathbb{R}^2.\label{eqvpws}
\end{eqnarray}
Here the sum extends up to $P$, the number of patches. The integrals along the contours are done counterclockwise. This is a family of equations depending on the parameter $\alpha$. In the limit $\alpha\rightarrow 0$ the equation approaches
the vortex patch problem 
of the 2D Euler equation whereas the limit $\alpha \to 1$ 
 approaches  the surface quasi-geostrophic equation.

 {The solutions of the $\alpha$-patches problem are numerically explored  in the range $0.5\leq \alpha \leq 1 $  by C\'ordoba {\it et al.} \cite{pnas} for the case of two patches and
evidence of  collapse is reported. 
Numerical  evidence of blow up is a  disputable issue, since 
 numerical calculations supporting this evidence do not constitute a formal proof.   This was the case,  for instance, for the 2D Euler equation \cite{Bu, DM} and still is for the 3D Euler equation \cite{kerr}. 
The current article  completes the results  by C\'ordoba {\it et al.}  describing in detail the
numerical method used to perform the simulations and providing benchmarks for this method.
Additionally the numerical results by C\'ordoba {\it et al.}  provide evidence that in the $\alpha$-patches problem, 
subsequent to a self-similar rescaling involving the parameter $\alpha$, the collapsing curves converge towards    a unique profile, thus  suggesting  the existence of a fixed point for the self-similar problem. While this problem was proposed in \cite{pnas}, it was  not thoroughly studied there. The current article is also focused in the study of collapse for $0.5\leq \alpha\leq 1$, particularly for the selection
$\alpha=0.7, 0.9$, but  makes progress by extending the numerical simulations to the problem reformulated in self-similar variables.  The advantage of using rescaled variables 
 is that the finite time collapse observed  by C\'ordoba {\it et al.}, becomes an asymptotic limit
in the pseudo-time $\tau$, and thus the 
evolution may be more accurately described. Our findings rule out the presence of a stationary graph   for the self-similar equation close to the convergent rescaled profiles.  However, a different but  exact fixed  point is described,  which is  found to play a role in separating collapsing from non-collapsing initial  data.}

This article is organized as follows:
Section 2 describes the problem; in particular the equations under study 
are introduced and the derivation of the self-similar equations is explained.  {General blow up  conditions are deduced for this setting.}
Section 3 reports on the numerical methods employed in the simulations. Section 4
provides an account of the results. Finally, the conclusions are presented in Section 5.

\section{The  equations}

In accordance with Zabusky et al. \cite{Z} the velocity  of a particle on the contour of a patch for the 2D vortex patch problem is obtained 
by inverting the relation between the streamfunction and the vorticity. For the $\alpha$-patches  problem
this relation is generalized  by 
$\theta= (-\Delta)^{1-\frac{\alpha}{2}}\psi $, where $\psi$ is the streamfunction and $\theta$ is the scalar that takes a constant value within the patch. In the limit 
$\alpha\rightarrow 0$, which corresponds to the 2D Euler equation, 
the scalar plays the role of vorticity and in the limit $\alpha\rightarrow 1$
(the quasi-geostrophic equation) the scalar  corresponds to the potential temperature.
For $0<\alpha<1$ one obtains that the velocity of a particle on the contour is given by:
\begin{equation}
u(\overrightarrow{x}(\gamma ,t),t)=\frac{\theta_{1}}{2\pi }    \ointctrclockwise_{C(t)}\frac{%
\frac{\partial \overrightarrow{x}}{\partial \gamma' }(\gamma ^{^{\prime }},t)%
}{|\overrightarrow{x}(\gamma ,t)-\overrightarrow{x}(\gamma ^{^{\prime
}},t)|^{\alpha }}d\gamma ^{^{\prime }}.\label{eqvp1}
\end{equation}
Here $\overrightarrow{x}(\gamma ,t)$ denotes  the position of a particle on the contour $C(t)$. This contour is 
parametrized by $\gamma $, and in our convention the integral along  it is done counterclockwise. Here $\theta_1={\theta^*}c_{\alpha}$, where 
the factor  
$\displaystyle{c_{\alpha}= \frac{\Gamma(\frac{\alpha}{2})}{2^{1-\alpha} \Gamma(\frac{2-\alpha}{2})}}$ 
results from inverting the operator $(-\Delta )^{1-\frac{\alpha }{2}}$ and
${\theta^*}$ is the value of scalar in the patch.  Changing the sign in the value of ${\theta^*}$ is equivalent to reversing the circulation along the integral in our convention.

The contour dynamic equation is eventually obtained by replacing  the velocity of a particle by the time derivative of its trajectory, {\it i.e.}:
\begin{equation}
\frac{d \overrightarrow{x}(\gamma ,t)}{d t}=u(\overrightarrow{x%
}(\gamma ,t),t)\text{ . }\nonumber
\end{equation}

This article focuses on the contour evolution of two patches.  The time evolution of a point $\vec{x}$ in any of the contours is  described by the equation:
\begin{eqnarray}
\frac{d\overrightarrow{x}(\gamma, t)}{dt}=\sum_{k=1}^{2} \frac{\theta_k}{2\pi}    \ointctrclockwise_{C_k(t)}
\frac{\frac{\partial \overrightarrow{x}_k}{\partial \gamma' }(\gamma ^{^{\prime }},t)  d\gamma' }{|\overrightarrow{x}(\gamma,t) -
\overrightarrow{x}_k(\gamma',t)|^{\alpha}}, \label{eqvpws2}
\end{eqnarray}
In the  limit  $\alpha =1$, which is the case of the quasigeostrophic equation,  local existence has been reported in \cite{Ro, paco}. In this case
one should use (see \cite{Ro}) the following formula for the
velocity (in the one-patch case):
\begin{equation}
u(\overrightarrow{x}(\gamma ,t),t)=\frac{\theta _{1}}{2\pi }   \ointctrclockwise_{C(t)}   \frac{%
\frac{\partial \overrightarrow{x}}{\partial \gamma' }(\gamma ^{^{\prime }},t)-
\frac{\partial \overrightarrow{x}}{\partial \gamma }(\gamma ^{^{ }},t)%
}{|\overrightarrow{x}(\gamma ,t)-\overrightarrow{x}(\gamma ^{^{\prime
}},t)|^{\alpha }}d\gamma ^{^{\prime }}\label{eqvp2}
\end{equation}
This equation eliminates the tangential component of the velocity, as only the normal component 
is able to deform the curve, thereby  avoiding  divergent integrals.

 \begin{figure}
\hspace{1cm}\resizebox{0.5\textwidth}{!}{\includegraphics{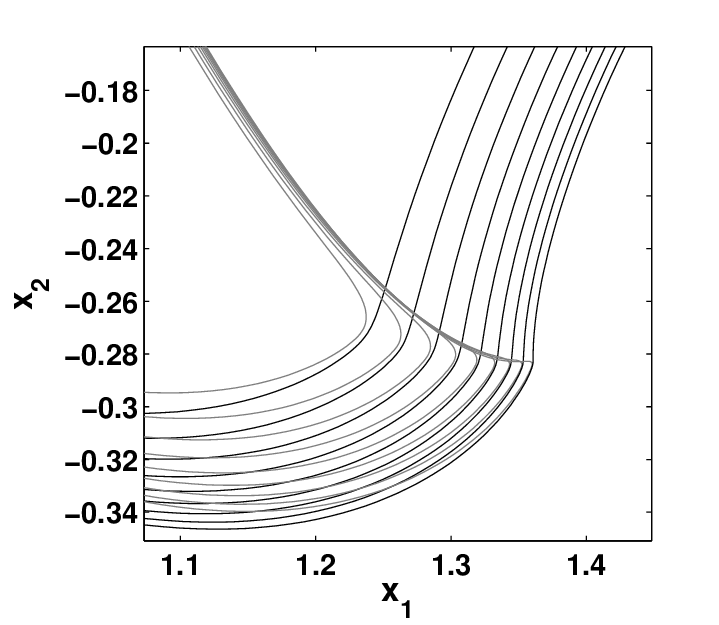}}\resizebox{0.5\textwidth}{!}{\includegraphics{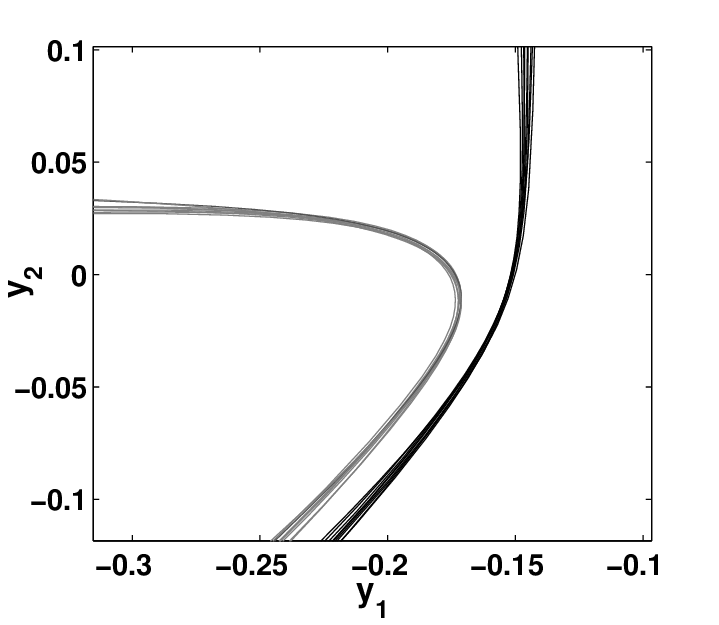}}
\caption{a) Zoom around the collapse point of two initial patches that evolve according to equation (\ref{eqvpws2}) for $\alpha=0.7$; b) the rescaled profiles near the collapse show the convergence towards a $\tau$ independent pattern.}
\label{fig:coll}
\end{figure}

Simulations  reported in \cite{pnas}  on Eq. (\ref{eqvpws2})  for $0.5\leq \alpha\leq 1$ show evidence of self-similar collapse.
Figure \ref{fig:coll} illustrates the collapse for $\alpha=0.7$ of two ellipses with semi-axis 1.1 (along the horizontal direction) and 1 (along the vertical coordinate) centered  
at positions $(-1.25, 0)$ and $(1.25, 0)$, on which $\theta_k=-1$ and the integrals are counterclockwise circulating. 
The singularity is point-like as shown in the magnification of Fig. \ref{fig:coll}(a) for a sequence near the corner.  { This  result is consistent
with the absence of splash singularity proven for this problem by  Gancedo and Strain \cite{paco2}, as it is found that  the curvature also blows up. }
 The coordinates of the collapse point are denoted as $\overrightarrow{x}_*(t_*)$, where
$t_*$ is the collapse time. Numerical results  reported in \cite{pnas}   indicate that  scaling laws 
exist near the blow up that describe the evolution of the maximum curvature as:
\begin{eqnarray}
\kappa \sim \frac{C}{(t_*-t)^\frac{1}{\alpha}} \quad {\rm as
}\quad t\rightarrow t_* \label{eq:k}
\end{eqnarray}
and also describe the minimum distance between contours as:
\begin{eqnarray}
d \sim {C}{(t_*-t)^\frac{1}{\alpha}} \quad {\rm as
}\quad t\rightarrow t_*\label{eq:d}
\end{eqnarray}
Fig. \ref{fig:coll}(b) shows  the profiles displayed in Fig. \ref{fig:coll}(a)  after rescaling   by a factor $1/(t_*-t)^{1/\alpha}$ and adjusting the output conveniently 
to achieve the coincidence in a $\tau$ independent pattern.

The result displayed in Fig. \ref{fig:coll}(b) led  C\'ordoba  et al. \cite{pnas} to propose an equation for  the rescaled  spatial variable $\overrightarrow{y}$: 
\begin{eqnarray}
\overrightarrow{x}(\gamma, t)-\overrightarrow{x}_*(t)&=&(t_*-t)^{\delta }(\overrightarrow{y}(\gamma, t)-\overrightarrow{y}_*),\label{eq:varxesc}
 \end{eqnarray}
where $\delta=1/\alpha$, and the velocity of the scaling function $\overrightarrow{x}_*(t)$ is given by
$$
\frac{d\overrightarrow{x}_*(t)}{dt}=\sum_{k}^{} \frac{\theta_k}{2\pi}  \ointctrclockwise_{C_k(t)}
\frac{\frac{\partial \overrightarrow{x_k}}{\partial \gamma' }(\gamma ^{^{\prime }},t)  d\gamma' }{|\overrightarrow{x}_*(t) -
\overrightarrow{x_k}(\gamma',t)|^{\alpha}}.
$$
The equation of motion for the rescaled space variable $\overrightarrow{y}$ is:
\begin{eqnarray}
(t_*-t)\frac{\partial \overrightarrow{y}}{\partial t }-\delta (\overrightarrow{y}-\overrightarrow{y}_*)=&&\nonumber\\
\sum_{k}\frac{\theta _{k}}{2\pi }    \ointctrclockwise_{\Upsilon_{k}(t)} \left( \frac{\frac{\partial \overrightarrow{y}_{k}}{\partial \gamma' }(\gamma ^{^{\prime }},t)}{|\overrightarrow{y}(\gamma ,t)-
\overrightarrow{y}_{k}(\gamma ^{^{\prime }},t)|^{\alpha}}-
\frac{\frac{\partial \overrightarrow{y}_{k}}{\partial \gamma' }(\gamma ^{^{\prime }},t)}{|\overrightarrow{y}_*-\overrightarrow{y}_{k}(\gamma ^{^{\prime }},t)|^{\alpha }} \right) d\gamma ^{^{\prime }}&&
\label{resc}
\end{eqnarray}
where  $\Upsilon_{k}$ are the rescaled contours. The new time variable 
\begin{eqnarray}
\tau &=&-\log(t_{*}-t),\label{eq:vartesc}
\end{eqnarray}
transforms the equation (\ref{resc})  into the self-similar equation:
\begin{eqnarray}
\frac{\partial \overrightarrow{y}}{\partial \tau }-\delta (\overrightarrow{y}-\overrightarrow{y}_*)=&&\nonumber\\
\sum_{k}\frac{\theta _{k}}{2\pi }    \ointctrclockwise_{\Upsilon_{k}(\tau)} \left( \frac{\frac{\partial \overrightarrow{y}_{k}}{%
\partial \gamma' }(\gamma ^{^{\prime }},\tau)}{|\overrightarrow{y}(\gamma ,\tau)-%
\overrightarrow{y}_{k}(\gamma ^{^{\prime }},\tau)|^{\alpha }}-
\frac{\frac{\partial \overrightarrow{y}_{k}}{%
\partial \gamma' }(\gamma ^{^{\prime }},\tau)}{|\overrightarrow{y}_*-\overrightarrow{y}_{k}(\gamma ^{^{\prime }},\tau)|^{\alpha }} \right) d\gamma ^{^{\prime }}&&
\label{resc1}
\end{eqnarray}
Like the non-rescaled equation, this equation admits a projection over $\overrightarrow{u}_n$,
the unitary normal component to the curve at $\overrightarrow{y}$,  which is useful for describing the $\alpha=1$ case,
\begin{eqnarray}
\frac{\partial \overrightarrow{y}_{\overrightarrow{u}_n}}{\partial \tau }=\left(\delta (\overrightarrow{y}-\overrightarrow{y}_*)\right)_{\overrightarrow{u}_n}+&&\nonumber\\
\left(\sum_{k}\frac{\theta _{k}}{2\pi } \ointctrclockwise_{\Upsilon_{k}(\tau)} \left( \frac{\frac{\partial \overrightarrow{y}_{k}}{%
\partial \gamma' }(\gamma ^{^{\prime }},\tau)}{|\overrightarrow{y}(\gamma ,\tau)-%
\overrightarrow{y}_{k}(\gamma ^{^{\prime }},\tau)|^{\alpha }}-
\frac{\frac{\partial \overrightarrow{y}_{k}}{%
\partial \gamma' }(\gamma ^{^{\prime }},\tau)}{|\overrightarrow{y}_*-\overrightarrow{y}_{k}(\gamma ^{^{\prime }},\tau)|^{\alpha }} \right) d\gamma ^{^{\prime }}\right)_{\overrightarrow{u}_n}&&
\label{resc2}
\end{eqnarray}
The projection evolves in time the same as the system (\ref{resc1}) ,  since  only the normal component deforms the curve. 
The  rescaled contours shown in Fig.  \ref{fig:coll}(b), which
seem to coincide over a unique curve, suggest that Eq. (\ref{resc2}) 
has a fixed point, also called a self-similar solution.  {Should an attracting fixed point linked to this pattern be found, it would be of great interest as
it  would provide a path toward rigorous proof  for the existence of blow up}. We will examine  this possibility further in the results section.

Interesting relations exist between the original and rescaled variables. 
For instance,  the blow up  in the original variables 
is transformed into an asymptotic behavior in the self-similar variables. 
In this way,  the collapse reported in \cite{pnas}, which occurs in a very small time interval in original variables, is mapped to an infinite 
interval in  rescaled variables, thus permitting a detailed monitoring of the blow up.
 We now consider the area of the patches;  this is related to the energy of the initial data and 
 is conserved in the original variables $(\overrightarrow{x}, t)$.
However, according to  Eq. (\ref{eq:varxesc}), it tends to grow in the rescaled variables $(\overrightarrow{y}, \tau)$. 
Other connections between the rescaled and non-rescaled variables
concern distance. The shortest distance, $d$, between
collapsing contours  (see Eq. (\ref{eq:d}) and Fig. \ref{fig:coll}(a)) tends to zero in the original variables, but  is not necessarily zero
in the new variables (see Fig. \ref{fig:coll}(b)). 
These results are easily justified. Since collapse is point-like,  the  distance between contours at the blow up time becomes zero just for  two  trajectories, one on each
contour, thus satisfying,
\begin{equation}
|\overrightarrow{x}_{C_1} (t)-\overrightarrow{x}_{C_2} (t)|\to 0, \,\,\,\,\,\, t\to t_* \label{eq:col2tr}
\end{equation}
The distance between these trajectories in self-similar variables is given by,
\begin{equation}
|\overrightarrow{x}_{C_1} (t)-\overrightarrow{x}_{C_2} (t)|=(t_*-t)^{\delta}|\overrightarrow{y}_{\Upsilon_1} (t)-\overrightarrow{y}_{\Upsilon_2} (t)| \label{eq:col2trrees}
\end{equation}
The factor $(t_*-t)^{\delta}$ in Eq. (\ref{eq:col2trrees}) confirms that the null distance 
between the original variables is satisfied even if trajectories on the rescaled
variables are at a finite distance above zero at the collapse time $t_*$. The above expression may be rewritten in terms of the pseudo-time $\tau$, as follows:
\begin{equation}
|\overrightarrow{x}_{C_1} (t)-\overrightarrow{x}_{C_2} (t)|=e^{-\tau \delta}|\overrightarrow{y}_{\Upsilon_1} (\tau)-\overrightarrow{y}_{\Upsilon_2} (\tau)| \label{eq:col2trreestau}
\end{equation}
As explained in the following section, the numerical technique used for simulations does not track individual particles on contours, but rather the  
contour as a whole, so the above collapse trajectories are not numerically integrated, 
only the contour that contains them. It is expected 
that near the collapse, the distance between these trajectories will be well represented by the evolution of the shortest distance
between  contours  which in \cite{pnas} is reported to evolve according to the expression:
\begin{equation}
|\overrightarrow{x}_{c_1} (t)-\overrightarrow{x}_{c_2} (t)|\sim  A \cdot (t_*-t)^{1/\alpha}, \,\,\,\,\,\, t\to t_* \label{eq:col3tr}
\end{equation}
Here $\overrightarrow{x}_{c_i} (t), i=1,2$ represents the set of points on the contours that are at a shortest distance as a function of time and $A$ is a constant.
Similarly,  the collapse may be also tracked in the self-similar variables  as the shortest distance between contours. 
This distance in  the self-similar variables is related to that in the non self-similar variables  by the expression:
\begin{equation}
|\overrightarrow{x}_{c_1}-\overrightarrow{x}_{c_2}|\sim  e^{-\tau \delta} |\overrightarrow{y}_{c_1}-\overrightarrow{y}_{c_2}|, \label{eq:col4tr}
\end{equation}
Let us define the minimum distance between contours in the  self-similar variables as:
\begin{equation}
 |\overrightarrow{y}_{c_1}-\overrightarrow{y}_{c_2}|\equiv D(\tau), \label{eq:col5tr}
\end{equation}
Eq. (\ref{eq:col4tr}) thus  indicates that a collapse, {\it i.e}, zero distance in the non self-similar variables, 
can be achieved  as long as the positive function $D(\tau)$ satisfies:
\begin{equation}
\lim_{\tau\to \infty}\frac{D(\tau)}{e^{\tau \delta} } \to 0, \,\,\,\,\,\,   \label{eq:col6tr}
\end{equation}
otherwise the distance between patches in the original coordinates could  not be collapsing.  In order to have a collapse condition, the positive function $ D(\tau)$ does not need to be asymptotically a constant; it can be time dependent or even a growing function, as far as the quotient in (\ref{eq:col6tr}) tends to zero.
This is consistent
with the collapse classification  in terms of self-similar variables reported  in \cite{eggersmarco}. 
It is shown in that work that  the asymptotic 
behavior of the collapsing data may be towards a fixed point (a constant $D(\tau)$), or towards $\tau$-dependent solutions
that may be either periodic or chaotic. In this article we  show that  seemingly collapsing data in the original variables have  an asymptotic regime
according to Eq. (\ref{resc1}),   which
shows no evidence of reaching a stationary 
regime.

Non-collapsing patches  always remain with 
finite curvature and at a finite distance  at any time in the original variables. In order for this to be the case,  the left hand side of Eq. 
(\ref{eq:col4tr}) needs to be finite, and this is only possible if 
\begin{equation}
\lim_{\tau\to \infty} D(\tau) \sim  e^{\tau \delta}. \,\,\,\,\,\,   \label{eq:colttr}
\end{equation}
If this asymptotic limit is not satisfied, as would be case for instance with a divergent limit   in Eq. (\ref{eq:col6tr}),  for instance, then the distance between non-collapsing patches in the non self-similar variables would diverge in finite time, and this 
is not a consistent outcome. It is verified in the Results section  that non-collapsing initial data asymptotically satisfy this condition, which in turn will be used  as a benchmark for the numerics.

\section{The numerical method}

 {The results discussed  by C\'ordoba et al \cite{pnas} on the numerical simulation of
Eq. (\ref{eqvpws2}) report evidence of blow up. Numerical  evidence of collapse is  always a disputable question, given that the
 numerical calculations supporting this evidence do not constitute a formal proof.   This was the case  for instance, for the 2D Euler equation \cite{Bu, DM} and still is for the 3D Euler equation \cite{kerr}. 
 The goal of this section is to provide full numerical details about the simulations performed in 
 \cite{pnas} for the $\alpha$-patches problem, as well as   extending   the simulations to the self-similar equation.  Benchmark examples supporting the correctness  of the results are discussed.

The time evolution of Eq. (\ref{eqvpws2}) is calculated by means of  contour dynamics.  This technique has been used in the past for the 2D Euler equations in several geophysical contexts \cite{ruso1,DM, DM2,ruso2,ruso3}.  
This method is particularly suitable because, contrary to   other numerical methods 
reporting singularity formation in the 2D Euler equations (see for instance \cite{Bu}), which turn out to be false \cite{Che, BeCo}, 
no singularity formation  is reported by authors using contour dynamics \cite{DM}. The consistency between theory and numerical
experiments supported by  this method in the 2D Euler problem ensures us of its robustness for the $\alpha$-patches problem. 
 However, the methodology discussed in \cite{DM2} cannot be straightforwardly applied  to the problem discussed  here because it deals with different equations. This section discusses   numerical details derived from handling the  $\alpha$-patches problem, which poses the difficulty of evaluating integrals that are more singular   than those  in  the vortex patches problem.  }

\subsection{The contour representation}

 {According to the methodology  described in \cite{DM2},   each contour $C_k$ is represented by a set of nodes $N_k$.
The curve between consecutive nodes  is interpolated by a cubic  spline}: 
\begin{equation}
\overrightarrow{x}_j (p) = \overrightarrow{x}_j + p \overrightarrow{t}_j + \eta_j (p)\overrightarrow{n}_j \label{eq:11}
\end{equation}
for $0 \leq  p \leq 1$ with $\overrightarrow{x}_j (0) = \overrightarrow{x}_j$ and $\overrightarrow{x}_j (1) = 
\overrightarrow{x}_{j+1}$, where:
\begin{eqnarray}
\overrightarrow{t}_j &=& (a_j , b_j ) = \overrightarrow{x}_{j+1} - \overrightarrow{x}_{j},  \,\,\, \overrightarrow{t}_j  \in \mathbb{R}^2\\
\overrightarrow{n}_j &=& (-b_j , a_j ),   \,\,\,  \overrightarrow{n}_j \in  \mathbb{R}^2\\
\eta_j (p) &=& \mu_j p + \beta_j p^2 + \gamma_j p^3, \,\,\,  \eta_j \in    \mathbb{R}.
\end{eqnarray}
The cubic interpolation coefficients $\mu_j$, $\beta_j$ and $\gamma_j$ are:
$$
 \mu_j = - \frac{1}{3} d_j \kappa_j - \frac{1}{6} d_j \kappa_{j+1}, \,\,\, \beta_j = \frac{1}{2} d_j \kappa_j, \,\,\, 
\gamma_j = \frac{1}{6}d_j (\kappa_{j+1} - \kappa_j ),
$$
where $d_j = |\overrightarrow{x}_{j+1} - \overrightarrow{x}_j |$ and
\begin{eqnarray}           
\kappa_j = 2 \frac{a_{j-1} b_j - b_{j-1} a_j}{|d^2_{j-1} \overrightarrow{t}_j + d^2_j \overrightarrow{t}_{j-1} |} \label{kappa}
\end{eqnarray}
is the local curvature defined by a circle through the three points, $x_{j-1}$, $x_{j}$, and $x_{j+1}$.
 The node  spacing in each contour is non-locally adjusted  at each time step depending on the curvature value.
Issues related to the density  of nodes in the curve will be addressed later. We now  explain 
 how  the above discretization transforms the system (\ref{eqvpws2}). The evolution of any point 
$\overrightarrow{x}(\gamma, t)$ on the contours is replaced by the evolution of a point $\overrightarrow{x}_j$ on the discrete curves, and
the integrals on the contour curves $C_k$ are replaced by the summation of integrals over the parameter $p$:
\begin{eqnarray}
\frac{d\overrightarrow{x}_j(t)}{dt}=\sum_{k=1}^{2} \frac{\theta_k}{2\pi} \sum_{i=1}^{N_k} \int_{0}^{1}
\frac{\frac{\partial \overrightarrow{x}_{i,k}(p,t)}{\partial p } \, dp }{|\overrightarrow{x}_j(t) -
\overrightarrow{x}_{i,k}(p,t)|^{\alpha}}, \label{eqvpws2d}
\end{eqnarray}
Here $\overrightarrow{x}_{i,k}(p,t)$ refers to a piece of curve computed as in Eq. (\ref{eq:11}) at a time $t$, and 
the additional subindex $k$ distinguishes the contour where the segment is placed.
 More abstractly, this expression may be written as:
\begin{eqnarray}
\frac{d\overrightarrow{x}_j(t)}{dt}=f_j(\overrightarrow{x}_j(t), C_1,C_2), \, \, \,\, \, \, f_j: \mathbb{R}^2 \to \mathbb{R}^2\label{compaut}
\end{eqnarray}
 This equation stands for any node $\overrightarrow{x}_j$ on any of the discretized contour curves.  So for the case of two contours on a plane,
 the system (\ref{eqvpws2d}) represents a set of $M=2 \times (N_1+N_2)$ coupled  ordinary differential equations, which
more compactly is rewritten as an autonomous system, as follows:
\begin{eqnarray}
\frac{d{\bf x}}{dt} = {\bf f} ({\bf x}), \,\,\,\,\,\,\,\, {\bf x} \in \mathbb{R}^M \label{aut}
\end{eqnarray}
The dimension $M$ of the system (\ref{aut}) is very large and typically non-constant, since the discretization does not force 
the number of points $N_k$ on the curve to be maintained for all times.
According to \cite{DM2}, the system (\ref{aut})
is integrated with an explicit  4$^{th}$ order Runge-Kutta method.  The time step $\Delta t$ in the Runge-Kutta method 
as suggested by C\'ordoba et al.  \cite{pnas} is chosen dynamically, since it has to be refined near the collapse time. As 
 reported in that work, the blow up concurs with
the formation of corners on the contours, and  the node spacing $\Delta x $ 
is reduced to represent corners properly. The time step   
is then adjusted as   $\Delta t= B \Delta x $, where $B$ is a constant  empirically tuned for different values of $\alpha$. 

\subsection{The evaluation of the contours integrals}

In the numerical simulation of the $\alpha$-patches problem, the most challenging part is the evaluation of 
each $f_j$ in Eq. (\ref{compaut}), 
as required by the Runge-Kutta method at each time step. In order to explain how this is achieved, we   next   focus on the contour integrals 
along each cubic contour segment (see Eq. (\ref{eqvpws2d})). In this  equation, the interpolation  (\ref{eq:11}) is replaced and
 integrals are left as follows:
\begin{flalign}
\int_{0}^{1}
\frac{ \left(   (\overrightarrow{t}_{i}+\mu_i \overrightarrow{n}_{i})+ (2 \beta_i p + 3 \gamma_i p^2 )\overrightarrow{n}_{i} \right) \, dp }{|
\overrightarrow{x}_i - \overrightarrow{x}_j + p \overrightarrow{t}_i + \eta_i (p)\overrightarrow{n}_i |^{\alpha}}&=&  \nonumber\\ (\overrightarrow{t}_{i}+\mu_i \overrightarrow{n}_{i})  \int_{0}^{1}
\frac{  \, dp }{|
\overrightarrow{x}_i - \overrightarrow{x}_j + p \overrightarrow{t}_i + \eta_i (p)\overrightarrow{n}_i |^{\alpha}}+&&\nonumber\\ 
 \overrightarrow{n}_{i} \int_{0}^{1}
\frac{   (2 \beta_i p + 3 \gamma_i p^2 ) \, dp }{|
\overrightarrow{x}_i - \overrightarrow{x}_j + p \overrightarrow{t}_i + \eta_i (p)\overrightarrow{n}_i |^{\alpha}}, &&\label{int}
\end{flalign}
The subindexes $k$ appearing in integrals of Eq. (\ref{eqvpws2d}) have been dropped to avoid cumbersome notation 
without loss of generality. 
As regards  the relative position of the variable 
$\overrightarrow{x}_j$ versus the fixed $\overrightarrow{x}_i$
 we classify the kind of integrals to be done as follows. 

{\bf 1.} The case $\overrightarrow{x}_j = \overrightarrow{x}_i$. In this situation, the first integral on the right hand side of Eq. (\ref{int}) has an integrand that becomes infinity at $p=0$. Nevertheless the principal value of the integral is defined if $\alpha<1$. In order to compute the PV, we rewrite the integral as follows:
\begin{eqnarray}
\int_{0}^{1}
\frac{     dp }{| p \overrightarrow{t}_i + \eta_i (p)\overrightarrow{n}_i |^{\alpha}}= \nonumber\\
\frac{1}{|\overrightarrow{t}_i|^{\alpha}} \frac{1}{(1+\mu_i^2)^{\alpha/2}}\int_{0}^{1}\frac{p^{-\alpha}dp}{\left(1+\frac{\beta_i^2p^2+\gamma_i^2 p^4+2\mu_i\beta_ip+2\mu_i\gamma_i p^2+2\beta_i \gamma_i p^3   }{(1+\mu_i^2)}\right)^{\alpha/2}}\nonumber
 \\=\frac{1}{|\overrightarrow{t}_i|^{\alpha}} \frac{1}{(1+\mu_i^2)^{\alpha/2}}\int_{0}^{1}
p^{-\alpha}\left(c_0+c_1p+c_2p^2\right)+\mathcal{O}(p^3)\nonumber\\
\sim \frac{1}{|\overrightarrow{t}_i|^{\alpha}} \frac{1}{(1+\mu_i^2)^{\alpha/2}}\sum_{n=0}^{10}\frac{c_n}{n-\alpha+1},
\label{intsing0}
\end{eqnarray}
Coefficients on the expansion may be easily computed  {with Maple or Mathematica or other algebraic manipulators}. In our numerical computations, we use an expansion up to the tenth order, which is very  
accurate as coefficients $c_n$ rapidly decay. We provide an explicit expression for $c_0$ and $c_1$.
$$
c_0=1,\,\,\, c_1=-\frac{\alpha \mu_i\beta_i}{1+\mu_i^2}
$$
These coefficients are also useful for approaching the second integral 
on the right hand side of Eq. (\ref{int}). In particular, this integral is approached by the summation:
\begin{eqnarray}
\frac{1}{|\overrightarrow{t}_i|^{\alpha}} \frac{1}{(1+\mu_i^2)^{\alpha/2}}\sum_{n=0}^{10}c_n \left(\frac{2\beta_i}{n-\alpha+2}+\frac{3\gamma_i}{n-\alpha+3}\right),
\end{eqnarray}

{\bf 2.} The case $\overrightarrow{x}_j = \overrightarrow{x}_{i+1}$. In this situation, the first integral on the right hand side of Eq. (\ref{int}) has an integrand that becomes infinity at $p=1$, as $x_i(p=1)=x_{i+1}$.
Results described in {\bf 1.} are  applicable here 
if integrals are rewritten for $p'=1-p$, 
$\tilde{\mu_i}=\mu_i + 2 \beta_i+3\gamma_i $, $\tilde{\beta_i}=- \beta_i-3\gamma_i $ and $\tilde{\gamma_i}=\gamma_i$.

{\bf 3.} The case $d_x=|\overrightarrow{x}_j -\overrightarrow{x}_{i}|>>0$. In this case the 
 first integral on the right hand side of Eq. (\ref{int}) 
may be approached by a series expansion as follows:
\begin{eqnarray}
\int_{0}^{1}
\frac{  \, dp }{|
\overrightarrow{x}_i - \overrightarrow{x}_j + p \overrightarrow{t}_i + \eta_i (p)\overrightarrow{n}_i |^{\alpha}}=&& \nonumber \\
\frac{1}{d_x^{\alpha}}\int_{0}^{1}\frac{dp}{\left(1+\frac{p^2 |\overrightarrow{t}_i|^2+\eta_i(p)^2|\overrightarrow{n}_i|^2+d_t p + d_n \eta_i(p)}{d_x^2}\right)^{\alpha/2}}&\sim& \frac{1}{d_x^{\alpha}}\sum_{n=0}^{10} g_n
\end{eqnarray}
where $d_n=-(\overrightarrow{x}_j-\overrightarrow{x}_i)\cdot \overrightarrow{n}$ and $d_t=-(\overrightarrow{x}_j-\overrightarrow{x}_i)\cdot \overrightarrow{t}$, and coefficients are computed 
with a software for algebraic manipulation. Coefficients $g_n$ rapidly decay if the distance, $d_x$, between
$\overrightarrow{x}_j$ and $\overrightarrow{x}_i$ is large enough. We provide an explicit expression for 
$g_0$ and $g_1$:
$$
g_0=1,\,\,\, g_1=-\frac{\alpha(d_t+d_n \mu_i)}{2 d_x^2}, ...
$$
 The second integral on the right hand side of Eq. (\ref{int}) may also  be expressed in terms of coefficients $g_n$.

\begin{figure}
\hspace{2cm}\resizebox{0.6\textwidth}{!}{\includegraphics{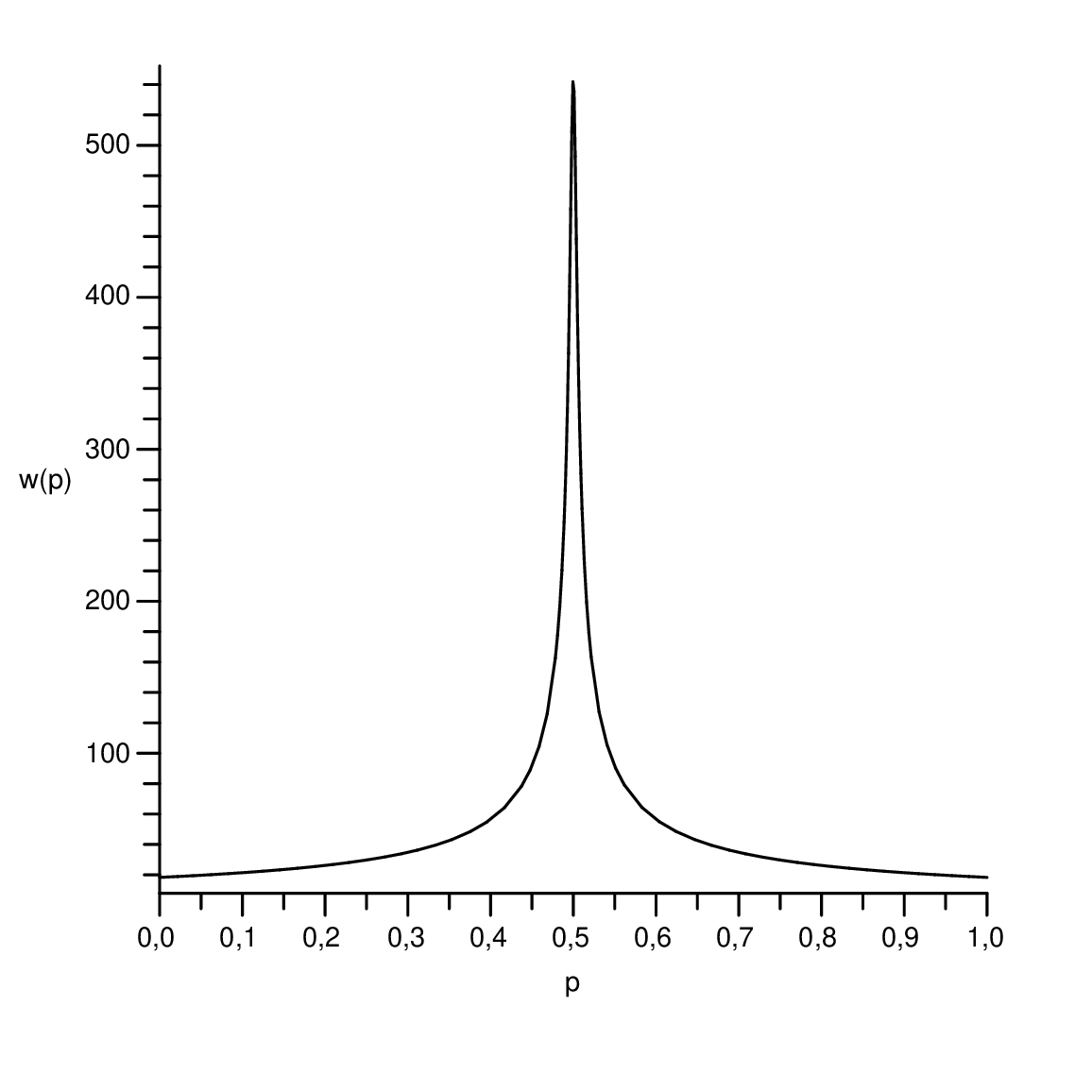}}\vspace{-0.5cm}
\caption{A graphic of the integrand in the first integral at the right hand side of  Eq. (\ref{int})
in the case $ |\vec{x}_j -\vec{x}_{i}| \sim 0$ it may vary steeply in the domain.}
\label{integrand}
\end{figure}

{\bf 4.} The case $ |\overrightarrow{x}_j -\overrightarrow{x}_{i}|\sim 0$. Here
the first integral on the right hand side of Eq. (\ref{int}) has an
 integrand $w(p)$ which may vary steeply in the domain taking for instance the appearance depicted in Fig. \ref{integrand}.
In this case, computing the area, $A$,  below the function:
$$
\int_0^1 w(p) dp =A
$$
with standard methods, such as a trapezoidal rule or Gaussian quadrature, may be inaccurate if the domain is not partitioned 
into a  grid that is fine enough to capture potential sharp features such as those displayed in  Fig. \ref{integrand}. On the other hand, these sharp features are not always present and if this is not the case 
a fine grid should not be necessary. In order to handle both possibilities accurately,  our perspective  for evaluating
 the above integral is to consider the following ordinary differential equation:
$$
\frac{dY}{dp}=w(p)
$$
with initial condition $Y(p=0)=0$. With this choice,  the integrated function will be such that $Y(p=1)=A$. The integration method 
used is a variable step  5$^{th}$ order Runge-Kutta method  (see \cite{nume}), thus ensuring it will appropriately track 
steep variations if required.
The evaluation of $w(p)$ must be conducted with great care to avoid round-off errors. This can be done correctly if the denominator
 in Eq. (\ref{int}) is computed as follows: First the polynomial in $p$ is expanded and the coefficient of 
each monomial in $p$ is evaluated; secondly, these terms are added up.  

In the numerical code, an appropriate test distinguishes which one of these procedures
needs to be followed to evaluate an  integral.  Special care needs to be taken to distinguish between cases 3 and 4. 
In particular,  choice 4 is adopted if the distance between $\vec{x}_j$ and $\vec{x}_i$  is below a  factor $f$ times the distance between
$\vec{x}_i$ and $\vec{x}_{i+1}$. The factor  typically ranges between 4 and 10 and is fixed by setting $f=1/\sqrt{Q}$, with $Q$ ranging from 0.05 to 0.01.

The evolution of the quasi-geostrophic equation (i.e the case  $\alpha=1$), as reported 
in Eq. (\ref{eqvp2}), is computed by projecting the vector field of Eq.  (\ref{eqvpws2})
over the unitary normal component to the curve at $\overrightarrow{x}_j $. This means that numerically 
the  first integral on the right hand side of Eq. (\ref{int}) does not need  to be computed (in fact, 
it is a divergent integral (see \cite{Ro})), since
it is a tangential component to the curve  at $\overrightarrow{x}_j $. The second integral  on the right hand side of Eq. (\ref{int}) needs to be projected over the unitary normal component to the curve, and this is easily done. The evolution of the contours
for $\alpha<1$ may be indistinctly computed either with Eq.  (\ref{eqvpws2}) or with Eq. (\ref{eqvp2}).

\subsection{A benchmark on the contour integrals}

The correctness and accuracy of the numerically computed contour integrals is verified by evaluating
$f_j(\overrightarrow{x}_j, C_1)$  for the case in which $\overrightarrow{x}_j=(1,0)$, $\theta_1=-1$ and there is only one contour, the circle of radius 1 centred at the origin.  This choice is appropriate because it can be  compared with the results obtained  from 
the natural parametrization in $\Theta$, i.e. $C_1=(\cos(\Theta), \sin(\Theta) ), \,\,\,\,\Theta \in [0, 2 \pi]$.
In this case, the integral to be evaluated is:
$$
f_j(\overrightarrow{x}_j, C_1)=\frac{-1}{2\pi}\int_0^{2\pi}\frac{(-\sin \Theta, \,\cos \Theta) \,\,\, d \Theta }{ \left(\sqrt{(1-\cos \Theta)^2+\sin^2 \Theta }\right)^{\alpha}}.
$$
which can be numerically evaluated with  Maple forcing 15 digits of precision, 
giving $f_j(\overrightarrow{x}_j, C_1)=( 9.549296\cdot 10^{-15},  -8.4000655 \cdot 10^{-1})$. The exact value of the first component in $f_j$ is zero, 
which is consistent with Maple forced precision.

 The  numerical evaluation explained in Section 3.2 has been checked for this one circle case with a number of points along the contour ranging from 100 to 450
 and tolerances $Q $ ranging from 0.05 to 0.01. Errors in the evaluation of $f_j$ extend from $\sim 10^{-5}$ (in the less favorable case with 100 points and $Q=0.05$)
 to   $\sim 5 \cdot 10^{-7}$  (in the more favorable case with 450 points and $Q=0.01$)

The precision of the numerical method thus depends on the number of points on the discretized curve, and also on the tolerance $Q$. 
In our results, we have typically used parameters ensuring accuracy of $f_j$   up to the 5$^{th}$ digit, which is  consistent
with the 4$^{th}$ order Runge-Kutta method employed to evolve {\bf x}, and with the precision expected from the redistribution procedure reported in the following subsection.

The numerical method is subjected to further analysis and tests, as discussed in the Results Section.

\subsection{The contour evolution}

Once ${\bf f} ({\bf x})$ in Eq. (\ref{aut}) is conveniently approached, we are ready to apply the 4$^{th}$ order Runge-Kutta method to
step forward ${\bf x}$.  At every time step  the nodes  are redistributed on the contour to guarantee  its optimal representation. This means that the numerical method will not track individual trajectories, but rather the whole contour. 
 {The relocation of points is performed  according to the methodology reported in \cite{DM2}, which is reviewed in the Appendix.}

\subsection{The self-similar problem}

The evolution of the self-similar equation (\ref{resc1}) may also be computed numerically, since its discrete version
may formally be rewritten as:
\begin{eqnarray}
\frac{d\overrightarrow{y}_j(\tau)}{d\tau}=F_j(\overrightarrow{y}_j, \Upsilon_1,\Upsilon_2), \, \, \,\, \, \, F_j: \mathbb{R}^2 \to \mathbb{R}^2\label{compaut2}
\end{eqnarray}
The numerical evaluation of $F_j$ at each step of the 4$^{th}$ order Runge-Kutta method
may be easily achieved with the use of the algorithm developed for the function $f_j$ since it is satisfied that:
$$
F_j(\overrightarrow{y}_j, \Upsilon_1,\Upsilon_2)=\delta \overrightarrow{y}_j+f_j(\overrightarrow{y}_j, \Upsilon_1,\Upsilon_2)-f_j(\overrightarrow{0}, \Upsilon_1,\Upsilon_2)
$$
where  the choice $\overrightarrow{y}_*=\overrightarrow{0}$ has been considered.

\section{Results}

 {The  method proposed in Section 3 is used for exploring numerically the solutions of  the self-similar  equation (\ref{resc1}).
We first   prove the existence of a stationary solution to the self-similar problem by providing its exact expression.
This solution  is shown to be valid in the range $0<\alpha\leq 1$, and is tried as a benchmark of the numerical calculations.
Numerical calculations require the selection of specific values of $\alpha$. We  systematically explore the results for
$\alpha$  equal to 0.7 and 0.9. As we have not found any essential differences, we report our findings for $\alpha$= 0.7.}

\subsection{An exact self-similar solution}

\begin{figure}
\hspace{2cm}\resizebox{0.7\textwidth}{!}{\includegraphics{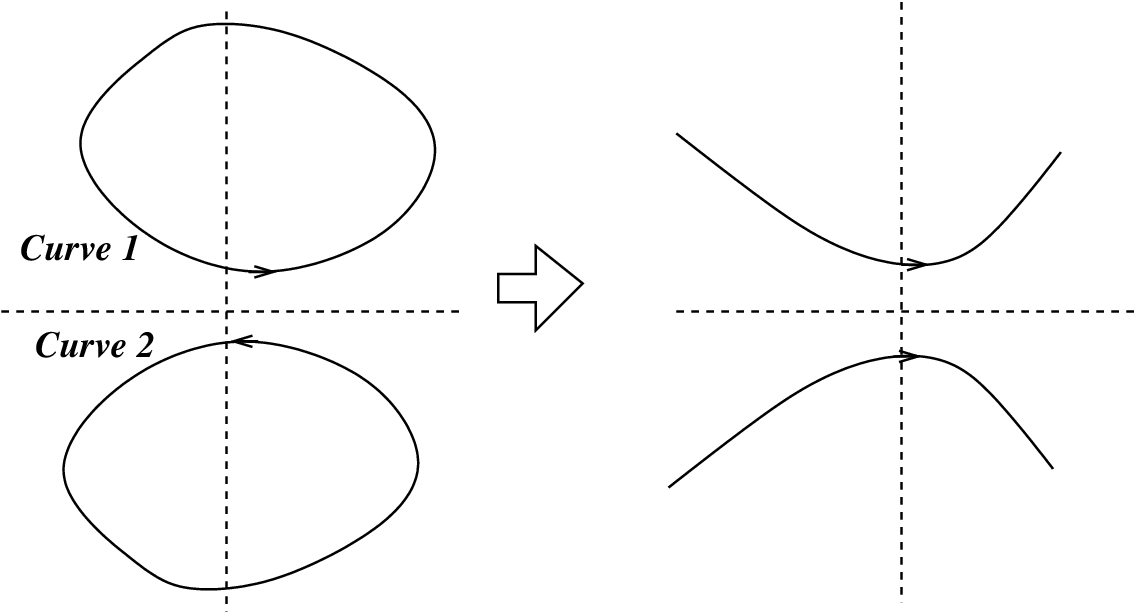}}\vspace{-0.5cm}
\caption{A scheme representing the parametrization of curves as functions.}
\label{funcion}
\end{figure}

  {At the collapse time $t_*$ the area within self-similar contours is infinite, as expected from rescaling (\ref{eq:varxesc}). 
 This suggests the possibility  that the stationary-like profiles of the rescaled curves take the form of a function $\Upsilon=\left(x,y(x)\right)$.}
When curves are parametrized as functions, Eq. (\ref{resc2}) is rewritten as:
\begin{eqnarray}
\frac{\partial \overrightarrow{y}}{\partial \tau } \rvert_{\overrightarrow{u}_n} = \delta (\overrightarrow{y}-\overrightarrow{y}_*) \rvert_{\overrightarrow{u}_n}   +\frac{\theta _{0}}{2\pi } \cdot &&\nonumber\\
\displaystyle ( \int_{\mathbb{R}}\frac{\frac{\partial \overrightarrow{y}_{1}}{%
\partial x' }(x',\tau) dx'}{|\overrightarrow{y}(x ,\tau)-
\overrightarrow{y}_{1}(x',\tau)|^{\alpha }}-
\int_{\mathbb{R}} \frac{\frac{\partial \overrightarrow{y}_{2}}{%
\partial x' }(x',\tau)dx'}{|\overrightarrow{y}(x ,\tau)-
\overrightarrow{y}_{2}(x',\tau)|^{\alpha }}-&&\nonumber\\
 \displaystyle \int_{\mathbb{R}}\frac{\frac{\partial \overrightarrow{y}_{1}}{\partial x' }(x,\tau)dx'}{|\overrightarrow{y}_*-\overrightarrow{y}_{1}(x',\tau)|^{\alpha }} }+
\int_{\mathbb{R}}\frac{\frac{\partial \overrightarrow{y}_{2}}{
\partial x' }(x',\tau)dx'}{|\overrightarrow{y}_*-\overrightarrow{y}_{2}(x',\tau)|^{\alpha }}    \displaystyle{)   \rvert_{\overrightarrow{u}_n}&&
\label{resc2r}
\end{eqnarray}
This expression takes into consideration that the counterclockwise integration direction on  the second curve  introduces   a minus sign when  it is replaced by an integration along  real numbers.  This is schematically represented
in figure \ref{funcion}.

\begin{figure}
\hspace{0.5cm}a)\resizebox{0.4\textwidth}{!}{\includegraphics{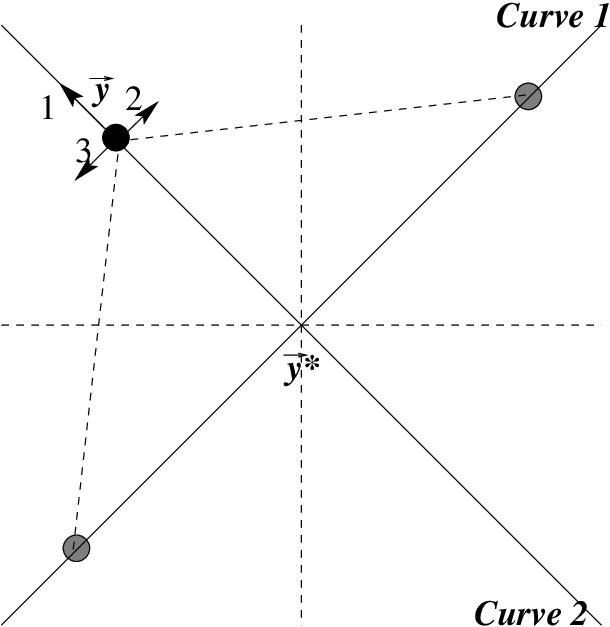}  } \hspace{1cm}b)\resizebox{0.55\textwidth}{!}{\includegraphics{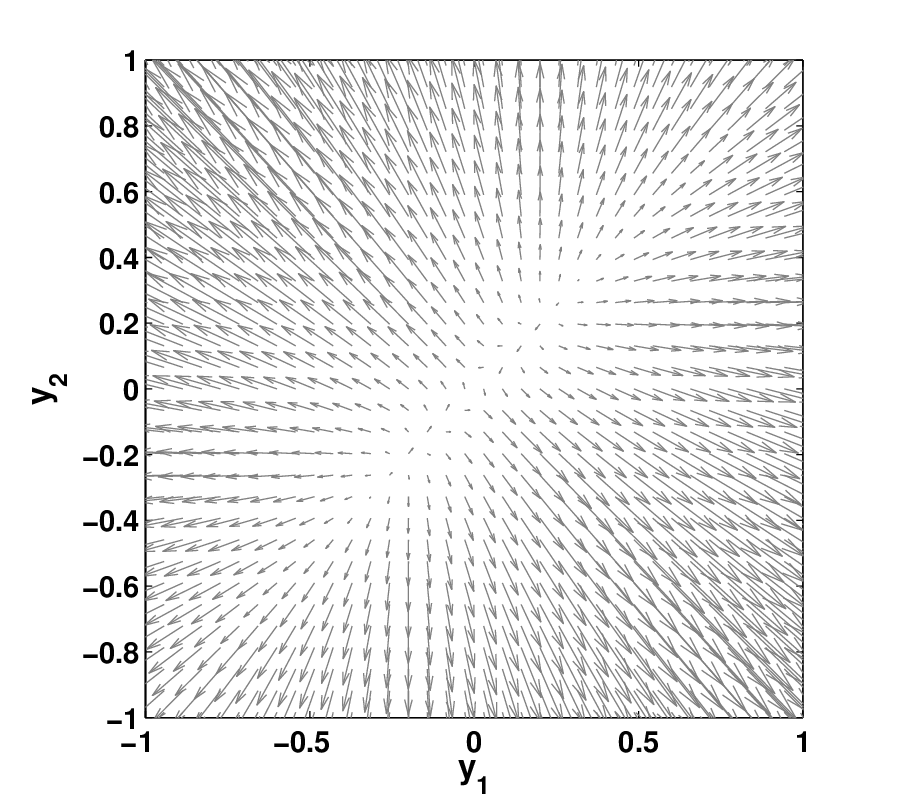}}
\caption{a) Sketch showing contributions to the point distinguished with a black  dot of different terms in Eq. (\ref{resc2r}).  
Each arrow is labelled with the number of the term it corresponds.  The geometry of the curves imposes cancellations of the normal components; b) velocity field numerically computed according to Section 3.2 which confirms the geometrical  argument.}
\label{geo}
\end{figure}

A fixed point
to Eq. (\ref{resc2r}) is:
\begin{eqnarray}
\Upsilon_1=(x, |x|),\,\,\,\,\,\,\,\,\,\Upsilon_2=(x, -|x|)\label{eq:selfsimsol}
\end{eqnarray}
These curves touch each other at one point, which is compatible with results on  regularity of the contours discussed in \cite{paco,paco2}. 
Figure \ref{geo}a) shows   how  geometrically over 
these particular curves the contributions of different terms of Eq. (\ref{resc2r}) to a point  $\overrightarrow{y}$,  indicated by a black dot, cancel each other out. The first term of Eq. (\ref{resc2r}) evaluated at the black dot, provides a null contribution on the normal direction, since before taking the normal component all the contribution of this term  is along the radial direction. This is highlighted with the arrow labelled with '1',  which is tangent to the curves expressed in  (\ref{eq:selfsimsol}). 
Arrows labelled  '2' and '3' show  the normal component of the small contributions of the second and third terms of    Eq. (\ref{resc2r}) integrated  along the small segments within the gray dots. These integrals depend only on the modulus of the distance between the gray and the black dots, and because of the symmetry they cancel each other out: they are equal but with different sign. 
The second and third terms  of    Eq. (\ref{resc2r}) have additional contributions to the black dot position $\overrightarrow{y}$, which do not cancel out but  which are tangent to the curve. Contributions from the  4$^{th}$ and 5$^{th}$ terms are the same for any potentially considered black dot $\overrightarrow{y}$, {\it i.e.} they do not depend on   $\overrightarrow{y}$, and they also cancel each other out.    {Clearly, the solution (\ref{eq:selfsimsol}) holds
for the whole range $0<\alpha\leq 1$, as our arguments do not depend on $\alpha$.}
The above argument is supported by the velocity field obtained from the numerical evaluation of the contour integrals   on the curves (\ref{eq:selfsimsol}), which is displayed in Figure \ref{geo}b). The calculation is performed
 according to the methodology explained in Section 3.2, and it serves  as a benchmark of the correctness of our numerics.

The fixed point in Eq. (\ref{eq:selfsimsol}) is not  unique, as
any rotation of it around the origin is also a stationary solution of Eq. (\ref{resc2}).  
Infinite stationary solutions exist that are related each other by this continuous transformation. 
In the context of infinite dimensional dynamical systems, there exist previous examples
 of families of fixed points  related by means of a 
group action (see for instance \cite{silvina}). It is  reported  in that work that this 
transformation implies   the existence of a null eigenvalue along a dynamically neutral eigenspace.

The evolution of  curves $\Upsilon_1$,  $\Upsilon_2$ is described in self-similar variables.
As  already pointed out, from Eq. (\ref{eq:varxesc}) it follows that in these coordinates
the area enclosed by the contours becomes infinite at the collapse time $t_*$, but  on the other hand blow up takes place when pseudo-time $\tau$ tends to infinity, as transformation (\ref{eq:vartesc}) confirms. This implies that
 at any finite stage  $\tau$ of the simulation of Eq. (\ref{resc1}),   the area of the patches  is finite, which makes 
 the numerical computation of its evolution  possible.

 \begin{figure}
a)\hspace{0.cm}\resizebox{0.55\textwidth}{!}{\includegraphics{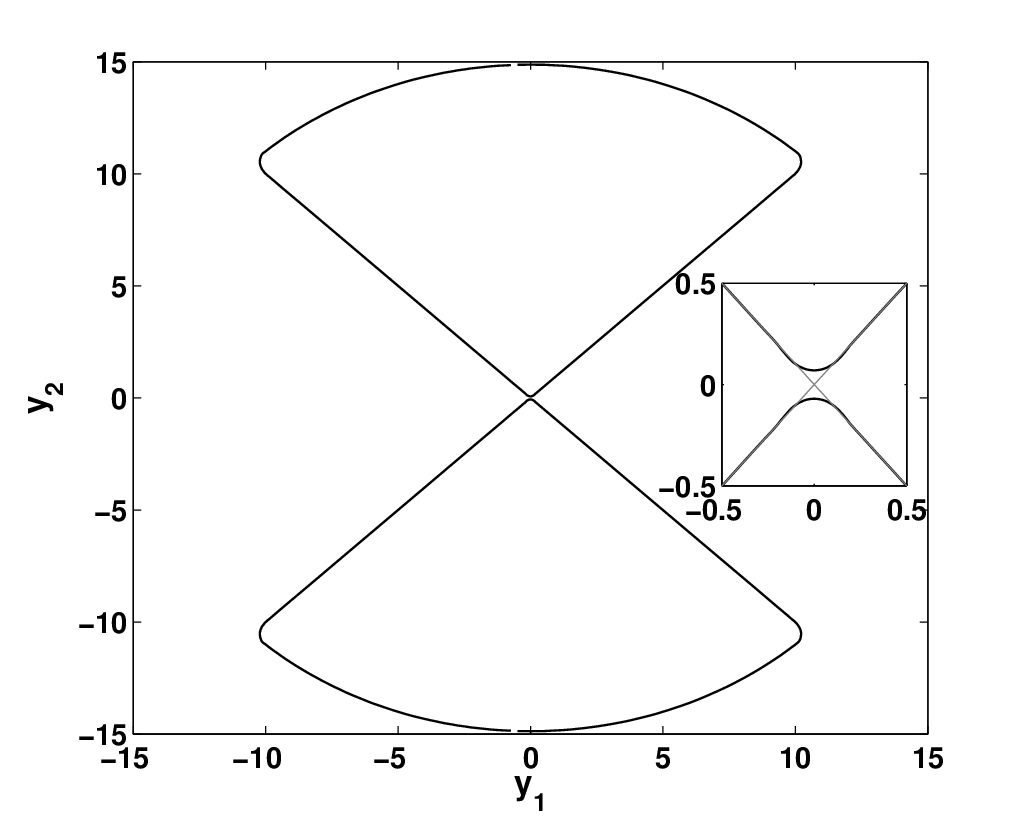}  } \hspace{0cm}b)\resizebox{0.55\textwidth}{!}{\includegraphics{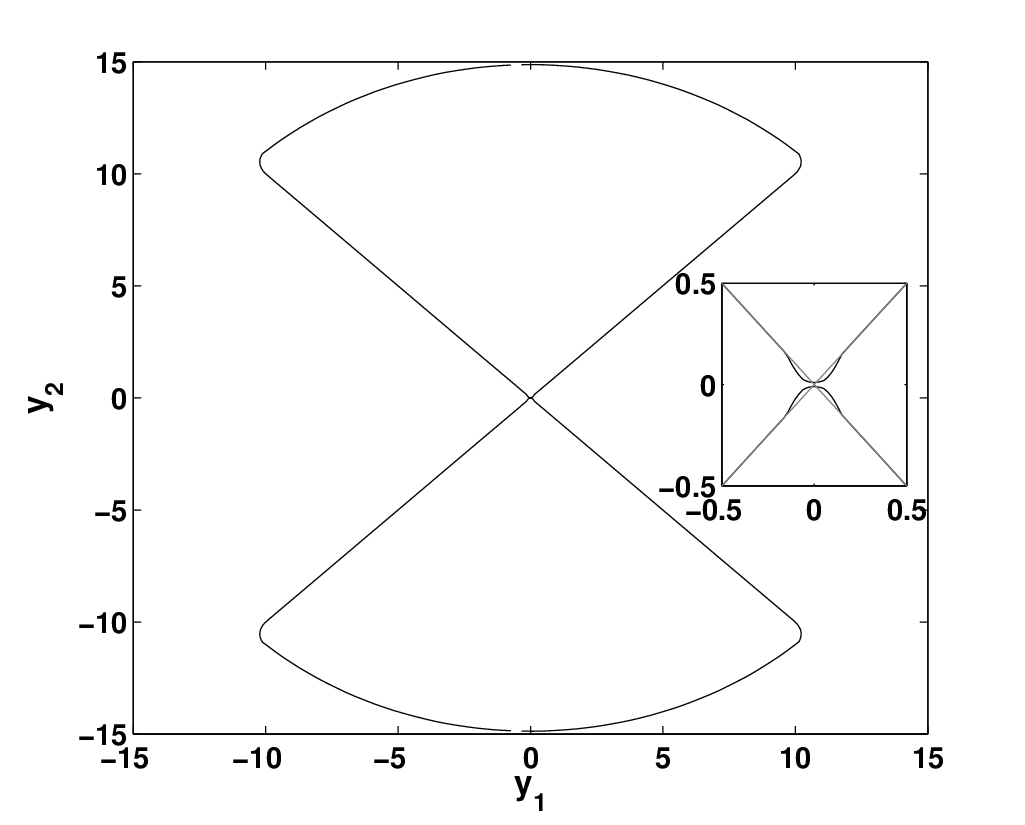}}
\caption{Two initial conditions  very close to the fixed point. The inset shows how are the initial patches close to the origin and gray lines stand for the fixed point.}
\label{fixper}
\end{figure}
We consider numerical simulations,  for the case $\alpha=0.7$, of curves very close to the fixed point in Eq. (\ref{eq:selfsimsol}) but  with slightly 
different perturbations near the origin.   { Figure  \ref{fixper} shows the initial data: the insets near the origin show
how  the perturbations  compare to the fixed point}. The closure of the curves far from the origin with partial circles does not affect the dynamics in the neighborhood of the origin.

\begin{figure}
a)\hspace{0.cm}\resizebox{0.55\textwidth}{!}{\includegraphics{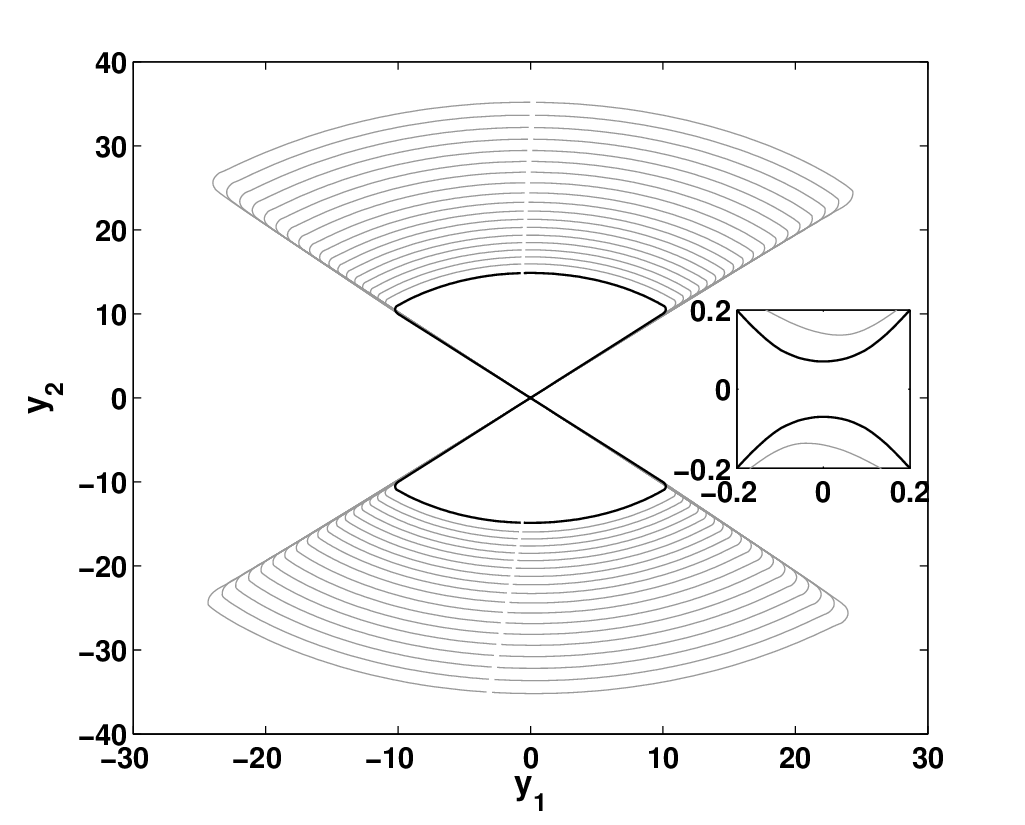}  } \hspace{0cm}b)\resizebox{0.55\textwidth}{!}{\includegraphics{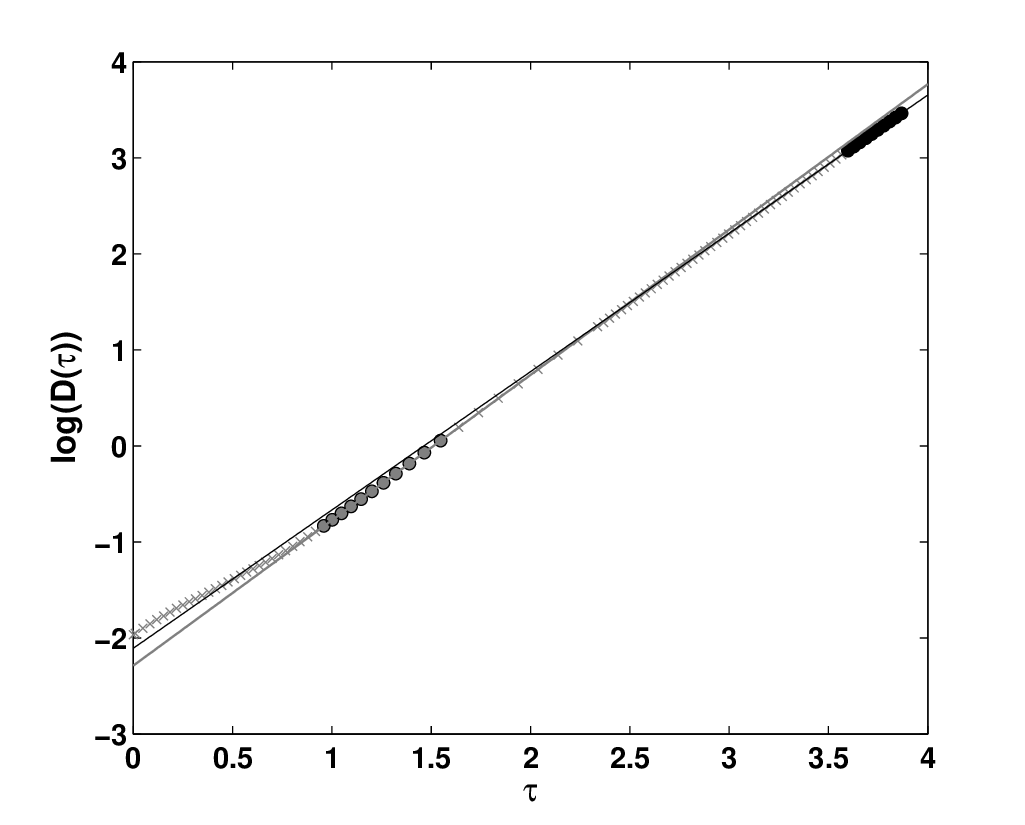}}
\caption{a) The evolution in  time of the initial  patches displayed in Figure \ref{fixper}a) ; b)  evolution of the logarithm of  the minimum distance between contours $D(\tau)$ versus the pseudo time $\tau$. The fitting line in gray 
is obtained for the gray dots and the fitting line black is obtained from the black dots.}
\label{fixpi}
\end{figure}

Our results on the evolution of the initial data in Figure  \ref{fixper}a) show evidence that the fixed point has repelling directions, in the sense that 
   there exist patches  very close to this point, as the ones in this panel,  that evolve in time, separating  from each other at a non-collapsing rate.
Figure \ref{fixpi}a) shows the  evolution in  time of the initial data represented in Figure  \ref{fixper}a). A zoom 
near the origin shows how the perturbed part evolves with respect to the fixed point. 
Figure \ref{fixpi}b) displays the time evolution of the logarithm of  the minimum distance between contours $D(\tau)$ versus the pseudo time $\tau$, thereby confirming a separation rate.
If this evolution is non-collapsing, according to Eq. (\ref{eq:colttr}), it  should asymptotically 
fit to a straight line with gradient $\delta=1/\alpha\sim 1.4286$. 
Two fittings have been performed over the data displayed in Figure  \ref{fixpi}b). One fitting for $\tau \in [0.9516, 1.54725]$  corresponds
 to the data marked with gray dots. The fitting line, shown in gray,  has a slope $s=1.5132$ while that obtained in the interval  $\tau \in [3.5982, 3.8682]$, displayed in black,  has a slope $s=1.4398$.  This confirms a decreasing gradient tendency for non-collapsing data which converges towards what is expected from Eq. (\ref{eq:colttr}).

\begin{figure}
a)\hspace{0.cm}\resizebox{0.55\textwidth}{!}{\includegraphics{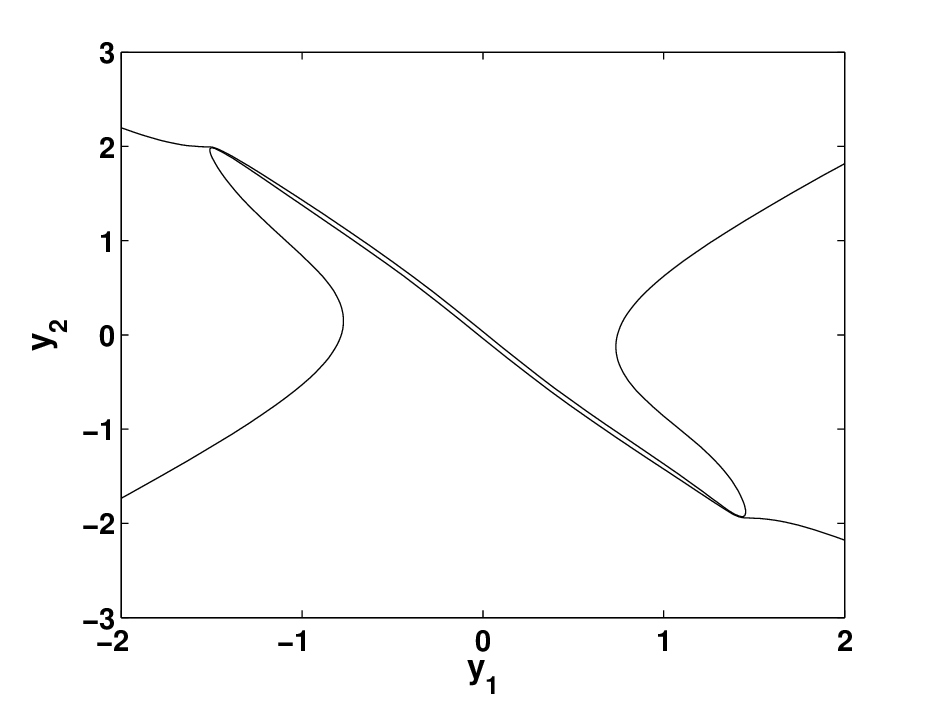}  } \hspace{0cm}b)\resizebox{0.55\textwidth}{!}{\includegraphics{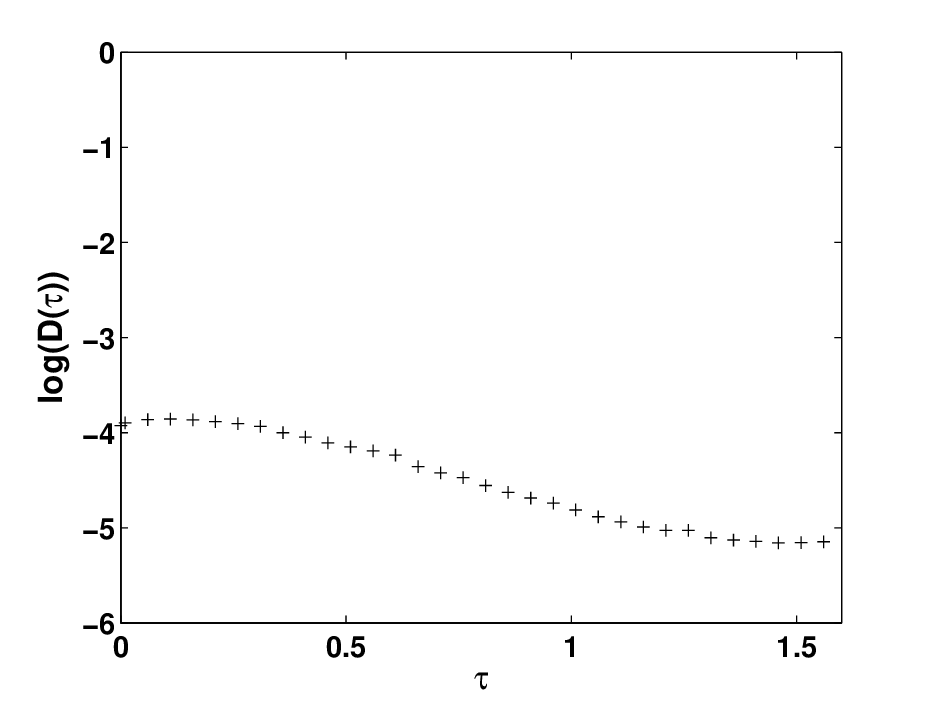}}
\caption{a)Two patches showing the formation of  two collapse points which move away along a  front;  b)  evolution of the logarithm of  the minimum distance between contours $\log(D(\tau))$ versus the pseudo time $\tau$.}
\label{distperiodic2}
\end{figure}

Next, perturbations to the fixed point (\ref{eq:selfsimsol}) are considered, such as that  shown 
in Figure \ref{fixper}b).  {The time evolution of this initial data is apparently  similar to  that observed  in Figure \ref{fixpi}a). However, in Figure \ref{distperiodic2}a) it is clearly distinguishable close to the origin of the formation of two collapsing points that move away  along a front.} The  evolution $\log(D(\tau))$ versus $\tau$ is displayed in Figure   \ref{distperiodic2}b). This signal indicates   the presence of a  collapse, in which $D(\tau)$ is kept finite.

\subsection{Non collapsing initial data}

\begin{figure}
\hspace{4.cm}\resizebox{0.55\textwidth}{!}{\includegraphics{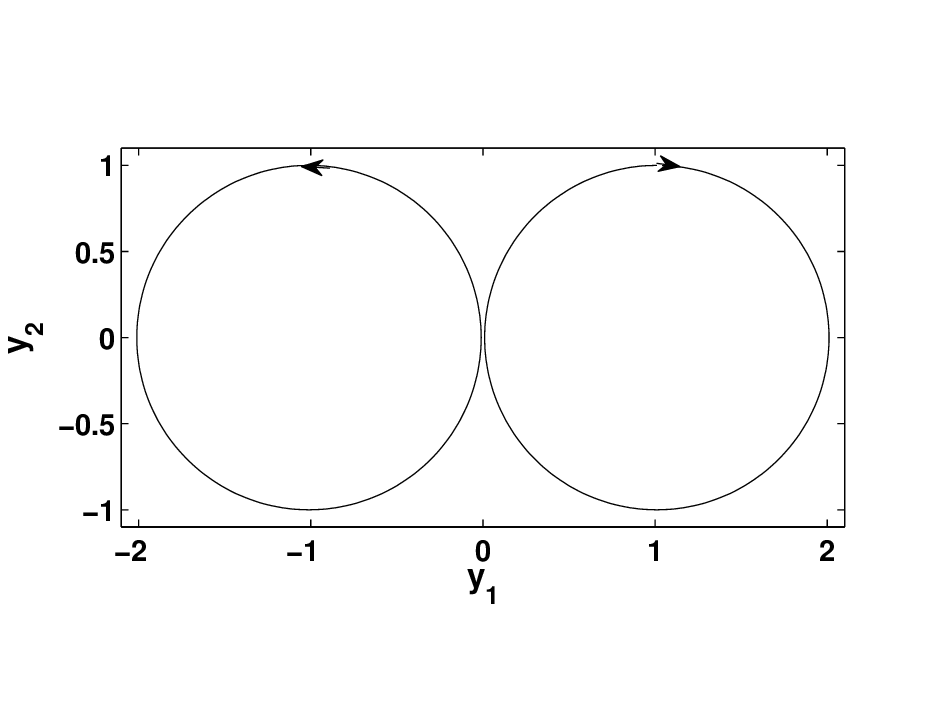}  }
\caption{Two initial patches of  counter-rotating circles.}
\label{noncounterrot}
\end{figure}

Further simulations are performed for data which are expected not to collapse, as for instance is the case of the two counter-rotating circles  displayed in Figure \ref{noncounterrot}. Here the circulation arrow refers to the integral circulation:  thus, for the 1st circle  this is done  as in equation (\ref{resc1}),  while for  the 2nd circle it is done in a way that is opposite to what is written in   there. This is equivalent to introducing a minus sign in front  of this term or considering $\theta_2=-\theta_1=1$.

\begin{figure}
\hspace{2.cm}\resizebox{0.7\textwidth}{!}{\includegraphics{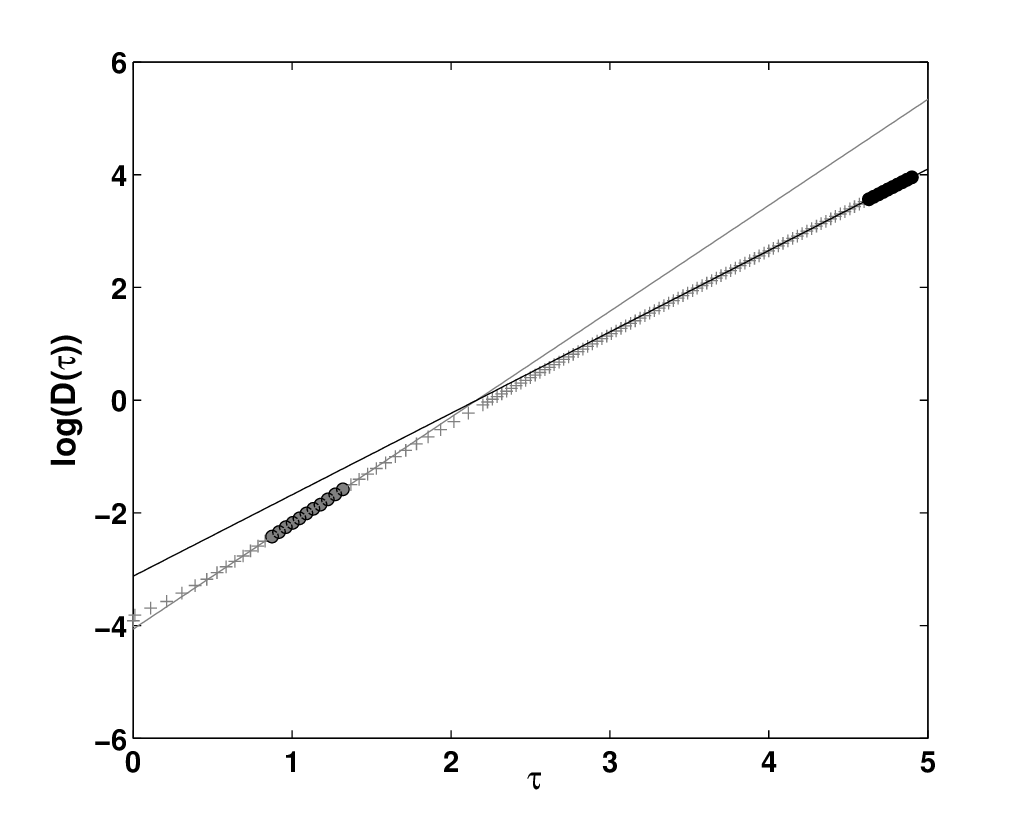}  }
\caption{Evolution of the logarithm of  the minimum distance between contours $D(\tau)$ versus the pseudo time $\tau$ for the counter-rotating circles. The fitting line in gray 
is obtained from the gray dots and the fitting line in black is obtained from the black dots.}
\label{circncrdt}
\end{figure}

The  pseudo time, $\tau$, evolution of the logarithm of  the minimum distance between contours $D(\tau)$  is shown in Figure \ref{circncrdt}. The  crosses in this graph indicate  distance  measurements from contours obtained in the simulation.  According to   Eq. (\ref{eq:colttr}), this non-collapsing initial data asymptotically should fit to a straight line with gradient $\delta=1/\alpha\sim 1.4286$. A fitting of the gray circles at early stages of the simulation 
in the interval $\tau \in [0.8736, 1.3193]$ supplies the slope $s=1.8793$, which indicates that the asymptotic regime has not been reached. A fitting at later times in the range $\tau \in [4.6297, 4.8997]$ 
confirms   a slope $s=1.4443$ decreasing towards the expected limit. These asymptotic results strengthen the consistency of the  numerics with exact findings, thus ensuring their accuracy.

\subsection{Collapsing initial data}

\begin{figure}
\hspace{0cm}a)\resizebox{0.55\textwidth}{!}{\includegraphics{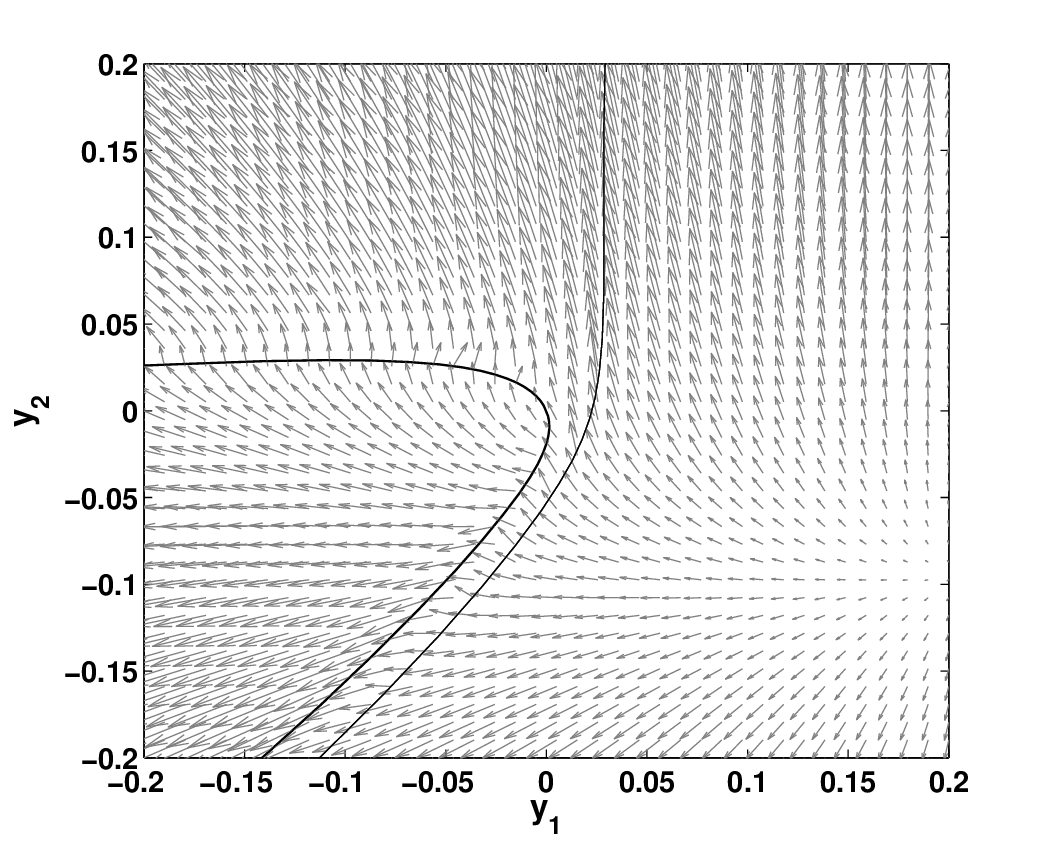}}b)\resizebox{0.55\textwidth}{!}{\includegraphics{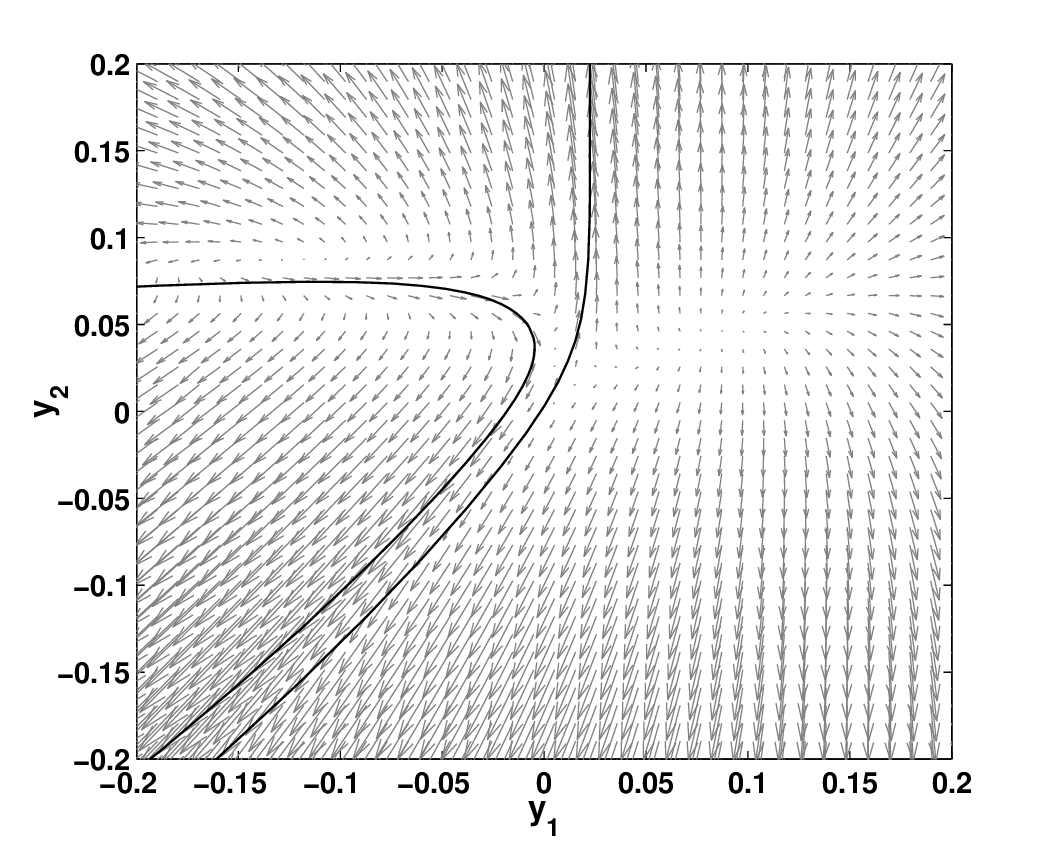}}
\caption{ Selfsimilar velocity field obtained near the collapse point.  a)  For ${\cal C}_2=(0,0)$; b) for ${\cal C}_2=(-0.0661, 0.0455)$.}
\label{qfp}
\end{figure}

\begin{figure}
\hspace{4.cm}\resizebox{0.6\textwidth}{!}{\includegraphics{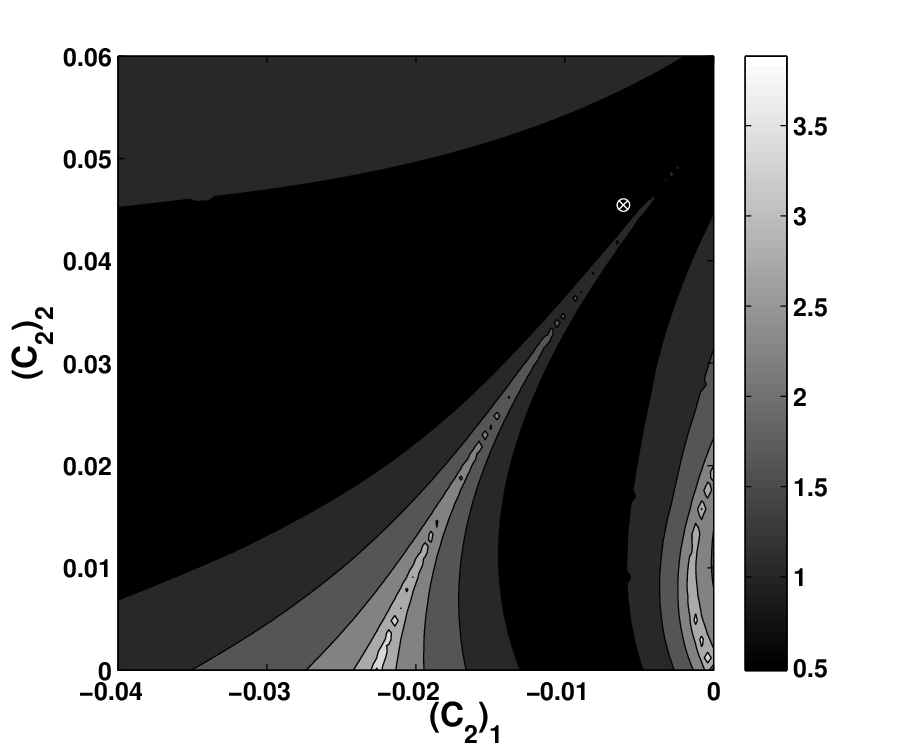}  }
\caption{Contours evaluating the sum  of the norm of the normal component of the velocity over selected points close to the corner on the curves $\Upsilon_1$ and $\Upsilon_2$ for different ${\cal C}_2$ values. The white circle highlights the position of the minimum at ${\cal C}_2=(-0.0661, 0.0455)$.}
\label{norm0i7}
\end{figure}

The simulations reported in Figure 1 summarize the findings reported in \cite{pnas} regarding the collapse of the $\alpha$-patches problem. This section discusses 
the evolution according to the self-similar equation (\ref{resc1}) of   initial data  -represented in Figure 1(a)- which is obtained  according to the non self-similar equation (\ref{eqvpws2}) 
at a time $t_{max}=6.81482306$, very close to the collapse time $t_*$.
According to the expression (\ref{eq:varxesc}),   an estimation of $t_*$, 
$\overrightarrow{x}_*(t_{max})$ and $\overrightarrow{y}_*$ is required to complete the transformation at  $t_{max}$. 
 In order to obtain the collapse time $t_*$,  a nonlinear fitting    is used in \cite{pnas}  which 
 closely follows the methodology reported in \cite{mancho}. This approach consists of  minimizing  with respect to all the parameters the quadratic differences between  the blowing-up curvatures obtained from the numerical simulation and those
evolving according to  the law (\ref{eq:k}). For the data reported in Figure 1, the estimated blow up time is  $t_*=6.8878$.  With regard  to  the evaluation of $\overrightarrow{x}_*(t_{max})$ and $\overrightarrow{y}_*$, both are unknowns: however, it is clear from 
the expression (\ref{eq:varxesc}) that once $t_*$ is estimated, both combine in a unique constant:
$$
{\cal{C}}=-\frac{\overrightarrow{x}_*(t_{max})}{(t_*-t_{max})^\delta}+\overrightarrow{y}_*
$$
which shifts the rescaled profiles:
\begin{eqnarray}
\overrightarrow{x}(\gamma, t_{max}) (t_*-t_{max})^{-\delta } ,\label{eq:rescprof}
 \end{eqnarray}
on the plane ${(y_1, y_2)}$. 
In practice  simulations of Eq.  (\ref{resc1}) consider that $ \overrightarrow{y}_*=\overrightarrow{0}$,  and for this reason fixing $\cal{C}$ at $t_{max}$,  fixes the unknown $\overrightarrow{x}_*(t_{max})$.
The constant $\cal{C}$ is split as follows: ${\cal C}={ \cal C}_1+{\cal C}_2$. Here ${\cal C}_1$ is  the constant 
that  shifts the point  with maximum curvature $\overrightarrow{x}^1_{max(\kappa)}$ in  curve 1 to zero, {\it i. e.},  ${\cal C}_1=-\overrightarrow{x}^1_{max(\kappa)}.$ 
Different ${\cal C}_2$ choices  produce significant changes in the velocity field near the origin.  
Figure \ref{qfp} confirms this point by showing the results for ${\cal C}_2=(0,0)$ and  ${\cal C}_2=(-0.0661, 0.0455)$.
In the second case, {\it i.e}, ${\cal C}_2=(-0.0661, 0.0455)$, one may observe how the rescaled curves behave as a quasi-fixed point of the velocity field. In other words,  close to the origin, 
the normal component of the  velocity field is almost zero on the curves. However, since this is something local,  they are not strictly speaking a fixed point to the equations.  
The choice ${\cal C}_2=(-0.0661, 0.0455)$  is such that for a selection of points on the curves $\Upsilon_1$ and $\Upsilon_2$ near the corner,  the constant minimizes the sum  of the norm of the normal component of the velocity over those points.
Figure \ref{norm0i7} shows how this sum varies for different values of the first and second components of ${\cal C}_2$. The white circle highlights the position of the minimizing ${\cal C}_2$. 

\begin{figure}
\hspace{4.cm}\resizebox{0.6\textwidth}{!}{\includegraphics{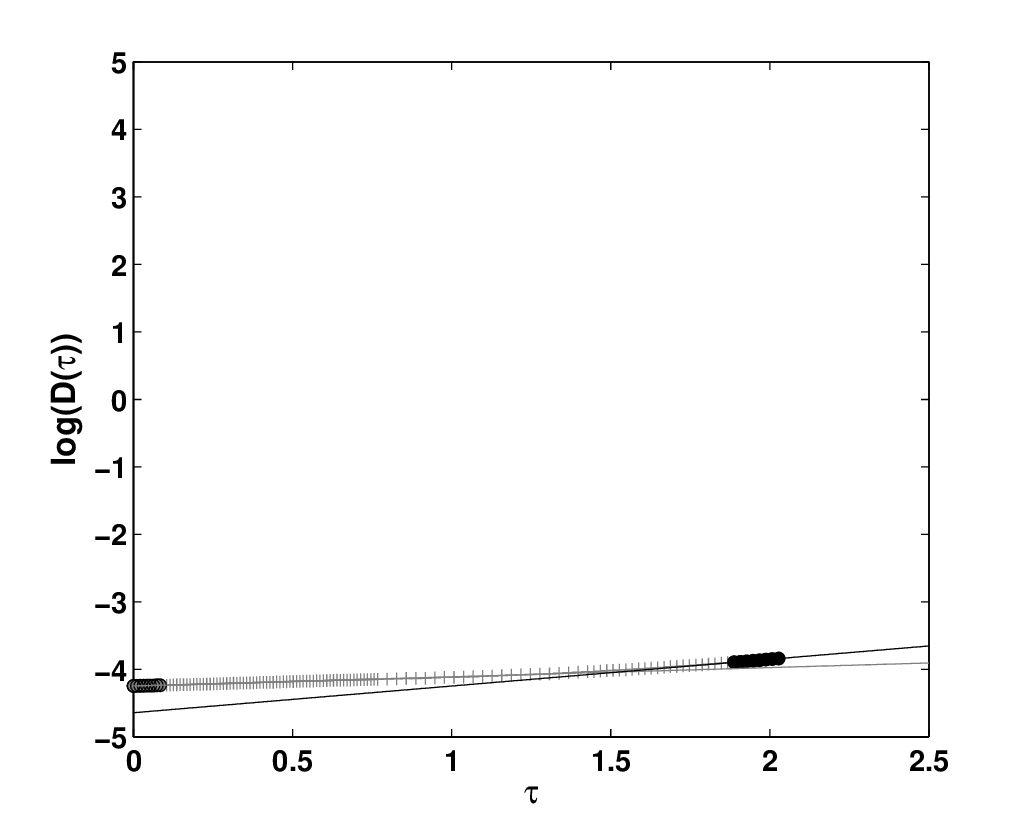}}
\caption{Evolution of the logarithm of  the minimum distance between contours $\log(D(\tau))$ versus the pseudo time $\tau$.}
\label{ev0i7qs}
\end{figure}

The numerical simulations summarized in Figure \ref{ev0i7qs} indicate that 
the time evolution of the curves  taking ${\cal C}_2=(-0.0661, 0.0455)$
 keep  the corner anchored  at the early stages,  when the curves remain at an almost stationary position, but eventually they separate and move apart.
The  crosses    in the graph  correspond  to the logarithm of  the minimum distance between contours $\log(D(\tau))$, obtained from the simulation.
  The  pseudo time, $\tau$, is displayed on the horizontal axis. A linear fitting of the gray circles 
in the interval $\tau \in [0,  0.0838]$ supplies the slope $s=0.1363 $ for the straight line in gray,   which confirms an almost stationary regime.  A fitting at later times, for the black bullets in the range $\tau \in [1.8875,  2.0275]$ 
confirms   the slope $s=0.3958$ for the straight line in black. This is well below the non-collapsing limit  $s=1/0.7=1.4286$, thus suggesting that this a collapsing solution, as reported in \cite{pnas}.

 {The evidence of  a quasi-stationary regime displayed throughout Figures  \ref{qfp}-\ref{ev0i7qs}   for the profiles shown in Figure \ref{fig:coll},
seem to rule out the possibility of finding   fully stationary graphs in the neighborhood of these profiles, as suggested in \cite{pnas}. The contribution of the  term $\delta (\overrightarrow{y})$ prevents any  cancelation (along the normal direction) on any possible elongation of these graphs -particularly far from the origin, when it becomes increasingly larger- with the exception of   the asymptotic branches  obeying the solution  (\ref{eq:selfsimsol}). But this fact brings us again to  Section 4.1 and the discussion therein about the manner in which   small perturbations to the exact self-similar solution evolve.}

\section{Conclusions}

The  $\alpha$-patches problem is studied mainly from the perspective of its re-formulation in self-similar variables, according to the approach proposed by
C\'ordoba et al. \cite{pnas}. 
 {Full details are given of the numerical implementation that  achieves a time simulation
of the problem, both in the original  and in the self-similar variables. Several benchmarks of the code are 
discussed throughout the paper. 

A main finding of this work is the proof of the existence of a family of stationary 
solutions to the self-similar problem by providing their exact expression. These solutions 
are shown to be valid in the range $0<\alpha \leq 1$. Simulations at selected $\alpha$ values ($\alpha=0.7$) of small perturbations 
around this solution indicate that it plays a role in separating collapsing from non-collapsing 
initial data. 

Numerical simulations on the self-similar variables allow a thorough analysis of the results reported   on \cite{pnas}.
 The 
advantage of using rescaled variables is that the finite time collapse becomes an asymptotic limit in the pseudo-time $\tau$, 
and thus the  evolution may be more accurately described. Near the blowup time,  rescaled contours are found in \cite{pnas} that seem to coincide over a unique graph.
 However, our results  show that  there is not such  
 stationary graph in the neighborhood of these convergent rescaled profiles. Instead, the simulations   indicate  that
they are related to the presence of a  quasi-stationary point, {\it i.e.},  they are  curves on which the normal component of the velocity is almost but not completely zero.
The separation rate between these curves is not constant, as expected in a simple collapse framework, but increases in the pseudo-time $\tau$. 
Nevertheless it is shown that, at least in the explored time ranges, the separation rate is slow enough to be consistent with blowup in the original variables.
Further simulations are discussed for collapsing and non-collapsing data.
The impossibility of performing  simulations for infinity pseudo-time means  that our results lie on extrapolations for observed tendencies, 
in spite of which they  still remain illustrative.}

\section*{Acknowledgements}
I wish to thank  A. C\'ordoba, D.  C\'ordoba, J. Eggers, C. Fefferman,  M. A. Fontelos  and J.L Rodrigo for 
their useful comments and suggestions. DC initiated me into $\alpha$-patches and suggested the problem discussed in this article. 
JE, CF and  MAF provided me
with valuable insights into the self-similar transformation. MAF helped me in the
evaluation of the integrals in Eq. (26) by means  of  series expansions. 
The computational part of this work was  performed on supercomputers at the BSC and  CESGA , and I am also grateful to these institutions for their valuable support.

This work is supported by Consolider I-MATH (C3-0104) and CSIC grants No. PI-200650I224 and OCEANTECH (No. PIF06-059) and MINECO grants: ICMAT Severo Ochoa project SEV-2011-0087 and  MTM2011-26696.

\section*{Appendix}

 {We review in  detail how the relocation of points is performed at each time step. 
According to \cite{DM2} a local density for each node is assumed:  }
$$
\rho_i=\frac{\hat{\kappa}_i}{1+\epsilon \hat{\kappa}_i/\sqrt{2}}\hat{\kappa}
$$
This form limits the minimum distance between nodes to nearly $\epsilon$ and the maximum computed curvature to $\epsilon^{-1}$.
We take $\epsilon=10^{-6}$.
Here $\hat{\kappa}_i=(\tilde{\kappa}_i+\tilde{\kappa}_{i+1})/2$ and
$$
\tilde{\kappa}_i=(\nu L)^{-1} (\breve{\kappa}_i L)^a+\sqrt{2}\breve{\kappa}_i
$$
 {We  consider $L=3$, $a=2/3$.  Here, $\nu$ is a parameter  that controls the accuracy of the representation. Typical values 
for $\nu$ in our study are 0.05, 0.03  or 0.01. Smaller values imply larger densities of points on the curve. For instance, the circle used for benchmark purposes in Subsection 3.3 with $\nu=0.01$ contains almost 450 points,
and with $\nu=0.05$ it contains around 100 points.}
The curvature $\breve{\kappa}$ incorporates non-local effects and is defined as:
$$
\breve{\kappa}_i=\sum_j \frac{d_j |\bar{\kappa}_j|}{h_{ij}^2}\left( \sum_j \frac{d_j}{h_{ij}^2}\right)^{-1}
$$
where $d_j=|\overrightarrow{x}_{j+1} -\overrightarrow{x}_{j}|$, $h_{ij}=|\overrightarrow{x}_{i}-(\overrightarrow{x}_{j+1} +\overrightarrow{x}_{j})/2|$ and $\bar{\kappa}_j=(\kappa_j+\kappa_{j+1})/2$ (see $\kappa_j$ defined in Eq. (\ref{kappa})) and the summation is over all nodes in the contour of the node $i$.

Given a  desired density $\rho_i$ between nodes $i$ and $i+1$,  it is possible to compute $\sigma_i=\rho_i d_i$, which is the fractional number of nodes to be placed between  $i$ and $i+1$. The node redistribution in a contour $k$ fixes one point on the
 curve, the first point, which will be kept fixed in the redistribution procedure. This implies that this point will 
evolve as a real trajectory of the system (\ref{aut}). 
The quantity $q$ is computed:
$$
q=\sum_{i=1}^{N_k}\sigma_i
$$
and let be $\tilde{N_k}$=$[q]+2$, the nearest integer to $q$ plus two. The $N_k-1$ 'old' nodes on the curve, are replaced by 
$\tilde{N_k}-1$ new nodes, which lie along the curved contour segments connecting the old nodes. 
Let $\sigma_i'=\sigma_i \tilde{N_k}/q$, so that $\Sigma_{i=1}^{N_k}\sigma_i'=\tilde{N_k}$. 
Then the positions of the new nodes $j=2...\tilde{N_k}$
are found successively by seeking $i$ and $p$ such that:
\begin{equation}
\sum_{l=1}^{i-1}\sigma_l'+\sigma_i'p=j-1
\end{equation}
and placing each new node $j$ between the old nodes $i$ and $i+1$ at the position $\overrightarrow{x}_i (p)$ given by Eq.
(\ref{eq:11}).

\end{document}